\documentclass{amsart}
\usepackage{amsmath}
\usepackage{amsthm}
\usepackage{amssymb}
\usepackage{amsfonts}
\usepackage{amsxtra}
\usepackage{epsfig}
\usepackage{verbatim}
\usepackage{enumerate}
\usepackage{amscd}

\theoremstyle{plain}
\newtheorem{thm}{Theorem}
\newtheorem{prop}[thm]{Proposition}
\newtheorem{cor}[thm]{Corollary}
\newtheorem{lem}[thm]{Lemma}

\theoremstyle{definition}
\newtheorem{defn}[thm]{Definition}

\theoremstyle{remark}
\newtheorem{remno}[thm]{Remark}
\newtheorem*{rem}{Remark}

\newcommand{\Real}{\mathbf R}
\newcommand{\C}{\mathbf C}
\newcommand{\N}{\mathbf N}
\newcommand{\Q}{\mathbf Q}
\newcommand{\Z}{\mathbf Z}
\newcommand{\A} {\mathbf{A}}
\newcommand{\h} {\mathbf{H}_3}
\newcommand{\oh} {\overline{\mathbf{H}}_3}

\begin{document}

\title[Denominators of Eisenstein
cohomology classes]{Denominators of Eisenstein cohomology classes
for ${\rm GL}_2$ over imaginary quadratic fields}

\author{Tobias Berger}

\address{Department of Pure Mathematics and Mathematical Statistics, Centre for Mathematical Sciences, University of Cambridge,
Cambridge CB3 0WB, United Kingdom}
\email{t.berger@dpmms.cam.ac.uk}



\begin{abstract}We study the arithmetic of Eisenstein cohomology classes
(in the sense of G. Harder) for symmetric spaces associated to
${\rm GL}_2$ over imaginary quadratic fields. We prove in many
cases a lower bound on their denominator in terms of a special
$L$-value of a Hecke character providing evidence for a conjecture
of Harder that the denominator is given by this $L$-value. We also
prove under some additional assumptions that the restriction of
the classes to the boundary of the Borel-Serre compactification of
the spaces is integral. Such classes are interesting for their use
in congruences with cuspidal classes to prove connections between
the special $L$-value and the size of the Selmer group of the
Hecke character.
\end{abstract}

\maketitle

\section{Introduction}
The relationship between the cohomology of an arithmetic subgroup
$\Gamma$ of a connected reductive algebraic group $G$ and the
automorphic spectrum of $\Gamma$ has been studied extensively. In
particular,  it is well-known that part of the cohomology can be
described by cuspidal automorphic forms. G. Harder initiated a
program to describe the entire cohomology in terms of cusp forms
and Eisenstein series (together with their residues and
derivatives). Using Selberg's and Langlands' theory of Eisenstein
series he constructed in \cite{HaGL2} a complement to the cuspidal
cohomology for the groups ${\rm GL}_2$ over number fields. These
Eisenstein classes can be described as cohomology classes with
nontrivial restriction to the boundary of the Borel-Serre
compactification of a symmetric space associated to $G$.

For arithmetic applications one would like to know if this
analytically defined decomposition
respects the canonical rational and integral structures on group
cohomology. Harder proved for ${\rm GL}_2$ that the decomposition
is, in fact,  rational. By the work of Franke and Schwermer \cite
{FS} a decomposition of the cohomology of a general reductive
group into cuspidal and Eisenstein parts and a rationality result
for the groups ${\rm GL}_n$ are now known. Harder also considered
the behavior with respect to the integral structure, in particular
the case when this decomposition is rational but not integral,
which corresponds to an Eisenstein class with integral restriction
to the boundary having a denominator. For a detailed exposition of
Harder's program we refer to \cite{HaSLN}.

We continue this analysis of the arithmetic of Eisenstein
cohomology classes in the case of ${\rm GL}_2$ over an imaginary
quadratic field $F$.  In this case, the associated symmetric space
is a 3-dimensional real manifold, and the cohomology in degrees 1
and 2 is the most interesting. We prove a lower bound on the
denominator of degree $1$ Eisenstein classes in terms of a special
$L$-value of a Hecke character, as conjectured by Harder. As an
example of the results proven, suppose $m \geq n \geq0$, let $p>
{\rm max} \{3,m\}$ be a prime split in $F$, and $\chi:
F^*\backslash \A_F^* \to \C^*$ a Hecke character with split
conductor coprime to $p$ of infinity type $z^{m+2} \overline
z^{-n}$ (see Theorem \ref{thm02} for the complete statement of our
result).  We construct an Eisenstein cohomology class ${\rm Eis}
\, \omega_{\chi}$ (for a coefficient system depending on $m$ and
$n$) that is an eigenvector for the Hecke operators at almost all
places such that the $p$-part of its denominator is divisible by
the $p$-part of $L^{\rm alg}(0,\chi)$. Here $L^{\rm alg}(0,\chi)$
is an integral normalization of the special $L$-value (see Theorem
\ref{Linteg}). In Proposition \ref{resinteg} we analyze when the
restriction of ${\rm Eis} \, \omega_{\chi}$ to the boundary of the
Borel-Serre compactification of the symmetric space is integral.
In particular, we prove this when $m=n$, $p>m+1$, and
$\chi^c(x):=\chi(\overline x)$ equals $\overline \chi(x)$ for all
$x \in \A_F^*$.

Such classes are interesting because of the implications for the
Selmer group of the $p$-adic Galois character associated to
$\chi^{-1}$: The situation here should be compared to the
classical Eisenstein series of weight 2 for $\Gamma_1(p)$ with a
character $\epsilon$ used by Ribet in \cite{Ri}. Its $q$-expansion
is $p$-integral and the constant term involves an $L$-value of
$\epsilon$. Via the congruence (of $q$-expansions) of the
Eisenstein series with a cuspidal Hecke eigenform Ribet proved the
converse to Herbrand's theorem. In our case the symmetric space is
not hermitian but one might try to use the integral structure
coming from Betti cohomology, as carried out for ${\rm GL}_{2/\Q}$
in \cite{HaPi} and \cite{CS}. If there exists an integral
cohomology class with the same restriction to the boundary as
${\rm Eis} \, \omega_{\chi}$ then our result shows that there
exists a congruence modulo $L^{\rm alg}(0,\chi)$ between ${\rm
Eis} \, \omega_{\chi}$, multiplied by its denominator, and a
cuspidal cohomology class. Via the Eichler-Shimura-Harder
isomorphism and the Galois representations attached to cuspidal
automorphic representations by the work of Taylor \textit{et al.}
(see \cite{T2}) one can then construct elements in the Selmer
group of $\chi^{-1}$ and obtain a lower bound on its size in terms
of $L^{\rm alg}(0,\chi)$. For this application of the results in
this paper in the case of constant coefficients see \cite{Be3}.

Note that for this application only the case $m=n$ is of interest
since cuspidal cohomology classes do not exist otherwise. Also,
since the interior cohomology for complex coefficients in degrees
1 and 2 are isomorphic, we restrict our study to degree 1. For an
analysis of denominators of degree 2 Eisenstein cohomology classes
associated to unramified characters  see \cite{F}.

We give a brief sketch of our proof of the lower bound on the
denominator in the special case of constant coefficient systems
(corresponding to $m=n=0$): In this case we can treat split or
inert primes $p>3$. Fix embeddings $F \hookrightarrow \overline \Q
\hookrightarrow \overline \Q_{p} \hookrightarrow \C$ and let
$\mathfrak{p}$ be the corresponding prime ideal of $F$ dividing
$p$. Let $G={\rm Res}_{F/\Q}({\rm GL}_{2/F})$ and $B$ the Borel
subgroup of upper-triangular matrices. For any (sufficiently
small) compact open subgroup $K_f \subset G(\A_f)$ let $S_{K_f}$
be the differentiable manifold $G(\Q) \backslash G(\A) / K_f
K_{\infty}$, where $K_{\infty}=U(2) \C^* \subset G(\Real)$. An
Eisenstein cocycle for $H^1(S_{K_f}, \C)$ is described by a pair
of Hecke characters $\phi=(\phi_1, \phi_2)$ with
$\phi_{1,\infty}(z)=z$ and $\phi_{2,\infty}(z)=z^{-1}$ and a
choice of a function $\Psi_{\phi_f}$ in the induced representation
$$V_{\phi_f, \C}^{K_f}=\{\Psi:G(\A_f) \to \C | \Psi(bg)=\phi_f(b) \Psi(g)
\forall b \in B(\A_f), \Psi(gk)=\Psi(g) \forall k \in K_f \}.$$ We
denote this Eisenstein cocycle by ${\rm Eis}(\Psi_{\phi_f})$. In
Section \ref{s3.2} we will make particular choices for
$\Psi_{\phi_f}$ (and corresponding $K_f$), the newvector
$\Psi^{\rm new}_{\phi_f}$ and the spherical vector
$\Psi^0_{\phi_f}$. We prove that $\Psi^{\rm twist}_{\phi_f}$, a
certain finite twisted sum of $\Psi^0_{\phi_f}$, is a multiple of
$\Psi^{\rm new}_{\phi_f}$ which will allow us to translate between
the two. The cohomology class $[{\rm Eis}(\Psi^0_{\phi_f})]$ is by
construction an eigenvector for the Hecke operators at almost all
places (see Lemma \ref{Hecke}) and we prove in Proposition
\ref{resinteg} that its restriction to the boundary is integral if
$\frac{L^{\rm alg}(-1, \phi_1/\phi_2)}{L^{\rm
alg}(0,\phi_1/\phi_2)}$ is, and proceed to show this is the case
if $(\phi_1/\phi_2)^c=\overline{\phi_1/\phi_2}$.

We know from the work of Harder that the cohomology class $[{\rm
Eis}(\Psi_{\phi_f})]$ is rational, i.e., it lies already in the
cohomology with coefficients in a finite extension of $F$. Since
we are interested in the $p$-adic properties we study, in fact,
its image in $H^1(S_{K_f}, \overline F_{\mathfrak{p}})$. The
denominator $\delta([{\rm Eis}(\Psi_{\phi_f})])$ of the Eisenstein
cohomology class  is the ideal by which it has to be multiplied to
lie inside the image of the cohomology with integral coefficients.
We prove that $$\delta([{\rm Eis}(\Psi^0_{\phi_f})]) \subseteq
(L^{\rm alg}(0, \phi_1/\phi_2)).
$$ By the functoriality of the evaluation pairing a cocycle represents an integral
cohomology class exactly when its pairing against all integral
cycles is integral. Explicit generators for the integral homology
are not known in our case, but we can obtain the desired lower
bound on the denominator by integrating ${\rm Eis}(\Psi_{\phi_f})$
against one carefully chosen integral cycle. The (relative) cycle
we use is motivated by the classical modular symbol: we integrate
along the path
$$\sigma:\mathbf{R}_{>0} \to {\rm GL}_2(\C)$$
$$ t \mapsto \begin{pmatrix} 1& 0 \\ 0 &t \end{pmatrix},$$
or rather a sum of such paths, one for each connected component of
$S_{K_f}$.

This ``toroidal" integral vanishes in general for
$\Psi^0_{\phi_f}$ but  we show that for $\Psi^{\rm
twist}_{\phi_f}$ the result, up to $p$-adic units, is
$$\int_{\sigma}{\rm Eis}(
\Psi^{\rm twist}_{\phi_f}) \sim \frac{L(0, \phi_1) L(0,
\phi_2^{-1})}{L(0,\phi_1/\phi_2)}.$$ We would like to conclude
from this that multiplication by at least $L^{\rm
alg}(0,\phi_1/\phi_2)$ is necessary to make our Eisenstein
cohomology class integral. For this we need to control the
$p$-adic properties of the numerator. To achieve this we use
results by Hida and Finis on the non-vanishing modulo $p$ of the
$L$-values $L^{\mathrm{alg}}(0, \theta \phi_i^{\pm 1})$ as
$\theta$ varies in an anticyclotomic $\Z_{q}$-extension for $q
\neq p$. We replace $\mathrm{Eis}(\Psi^{\rm twist}_{\phi_f})$ by
another ``twisted" version $\mathrm{Eis}^{\theta} (\Psi^{\rm
twist} _{\phi_f})$ for a finite order character $\theta$ of
conductor $q^r$, defined by
$$\mathrm{Eis}^{\theta} (\Psi^{\rm twist} _{\phi_f})(g)= \sum_{x \in
(\mathcal{O}_q/q^r)^*} \theta^{-1}(x) {\rm Eis}(\Psi^{\rm twist}
_{\phi_f})(g \begin{pmatrix} 1 & -\frac{x}{q^r} \\0 &1
\end{pmatrix}_q),$$ where $\mathcal{O}_q$ is the ring of integers of the completion of $F$ at $q$.
 The sum of paths making up the cycle is also weighted by values of
$\theta$. See Section \ref{3.2.1} for the definition of this cycle
$\sigma_{\theta}$. Up to units the result of this toroidal
integral is
$$\int_{\sigma_{\theta}} \mathrm{Eis}^{\theta}(\Psi^{\rm twist}_{\phi_f}) \sim \frac{L(0,
\phi_1 \theta) L(0, \phi_2^{-1}
\theta^{-1})}{L(0,\phi_1/\phi_2)}.$$ The results of Hida and Finis
allow us (under certain conditions on the conductors of the
$\phi_i$) to find  a character $\theta$ such that the numerator is
a $p$-adic unit. Apart from differences in the conditions on the
conductors Hida deals only with split $p$, whilst Finis also
treats
 inert $p$ for constant coefficients. Given a character $\chi$ satisfying certain assumptions
 we prove in Theorem \ref{thm02} the existence of characters
 $\phi_1$ and $\phi_2$ with $\chi=\phi_1/\phi_2$ for which the
 $L$-values in the numerator can be simultaneously controlled.
 This involves the construction of characters with prescribed
 ramification and a careful analysis of Artin roots numbers.

The twisting by $\theta$ also has the effect of making
$\mathrm{Eis}^{\theta}(\Psi^{\rm twist}_{\phi_f})$ vanish at the
$0$- and $\infty$-cusps of each connected component. By a result
of Borel (see Proposition \ref{Borel}) it therefore represents a
relative cohomology class with respect to these boundary
components. We prove that this relative cohomology class is again
rational and that its denominator bounds that of
$\mathrm{Eis}^{\theta} (\Psi^{\rm twist} _{\phi_f})$ from below.
We can therefore interpret the toroidal integral as an evaluation
pairing between relative cohomology and homology and deduce that
the ideal generated by $L^{\rm alg}(0, \phi_1/\phi_2)$ gives a
lower bound on the denominator of the relative cohomology class
represented by $\mathrm{Eis}^{\theta}(\Psi^{\rm twist}
_{\phi_f})$. We conclude the desired bound on the denominator of
$[{\rm Eis}(\Psi^0_{\phi_f})]$ by using the divisibilities $$
\delta([{\rm Eis}(\Psi^0_{\phi_f})]) \subseteq \delta(
[\mathrm{Eis}(\Psi^{\rm twist}_{\phi_f})]) \subseteq
\delta([\mathrm{Eis}^{\theta}(\Psi^{\rm twist} _{\phi_f})])
\subseteq \delta([\mathrm{Eis}^{\theta}(\Psi^{\rm twist}
_{\phi_f})]_{\rm rel}).$$

Our results generalize and extend the work in \cite{Ko} for
$F=\Q(i)$ and unramified $\phi_1/\phi_2$, where the toroidal
integral is calculated for the spherical vector. K\"onig proceeds
to show in his case that the $L$-value gives an upper bound on the
denominator.  Before this, Eisenstein cohomology for imaginary quadratic 
fields had been studied in \cite{Ha79}, \cite{Ha82}, and \cite{Wes}. 
Previous work on calculating or bounding
denominators for ${\rm GL}_2$ over $\Q$ and totally real fields
include \cite{chapsix}, \cite{Ka}, \cite{Maennel}, \cite{Mahn2},
\cite{CS}, and \cite{Wa}. \cite{Ka} and \cite{CS} also use
twisting techniques and a result by Washington on the
non-vanishing modulo $p$ of Dirichlet $L$-values in cyclotomic
towers. New about our method for getting a lower bound is that we
introduce the auxiliary cocycle $\mathrm{Eis}^{\theta} (\Psi^{\rm
twist} _{\phi_f})$ and prove that it represents a relative
cohomology class, which allows us to work just with the toroidal
integral, making the calculation of additional boundary integrals
as in \cite{chapsix}, \cite{Ka} unnecessary. Our method does not
allow to prove upper bounds because of the transition to the
finite twisted sum, but one might be able to get an upper bound by
applying this idea to prove a lower bound on the denominator of
the dual cohomology class in degree 2. In principle, our method
should extend to general CM-fields, where Hida's result is still
applicable. Since the arithmetically interesting classes appear in
the middle degrees this would, however, be notationally more
cumbersome (but see \cite{Maennel}).

These results generalize part of my thesis \cite{Be} under C.
Skinner at the University of Michigan, where this problem was
considered in the case of constant coefficient systems and split
$p$. The author would like to thank Thanasis Bouganis, Vladimir
Dokchitser, G\"unter Harder, Joachim Schwermer, and Chris Skinner
for helpful discussions and an anonymous referee for improvements
to the introduction and corrections in the statement of Theorem
\ref{Linteg}. This article was written during visits to the Max
Planck Institute in Bonn and the Erwin Schr\"odinger Institute in
Vienna. The author would like to thank both for their hospitality
and support.

\section{Notation and Definitions}
\subsection{Basic notation}
 Let $F$ be an imaginary
quadratic field, $\sigma$ its nontrivial automorphism,
$\mathcal{D}$ the different of $F$, and $d_F={\rm
Nm}(\mathcal{D})$ the absolute discriminant. For a place $v$ of
$F$ let $F_v$ be the completion of $F$ at $v$. We write
$\mathcal{O}$ for the ring of integers of $F$, $\mathcal{O}_v$ for
the closure of $\mathcal{O}$ in $F_v$, $\mathfrak{P}_v$ for the
maximal ideal of $\mathcal{O}_v$, $\pi_v$ for a uniformizer of
$F_v$, and $\hat{\mathcal{O}}$ for $\prod_{v \, {\rm finite}}
\mathcal{O}_v$. Complex conjugation is denoted by $z \mapsto
\overline z$. We use the notations $\A, \A_f$ and $\A_F, \A_{F,f}$
for the adeles and finite adeles of $\Q$ and $F$, respectively,
and write $\A^*$ and $\A_F^*$ for the group of ideles. Let $p>3$
be a prime of $\Z$ that does not ramify in $F$. Fix embeddings $F
\hookrightarrow \overline \Q \hookrightarrow \overline \Q_{p}
\hookrightarrow \C$ and let $\mathfrak{p}$ be the corresponding
prime ideal of $F$ over $p$.

\subsection{The algebraic group and symmetric spaces} \label{s2.2}
For any algebraic group $H/\Q$ and any ring $A$ containing $\Q$ we
write $H(A)$ for the group of $A$-valued points. We shall
abbreviate $H_{\infty}=H(\Real)$. We consider the algebraic group
$$G:=\mathrm{Res}_{F/\Q} ({\rm GL}_{2/F}).$$ The group $G_0/F={\rm
GL}_{2/F}$ contains the Borel subgroup of upper triangular
matrices $B_0$, its unipotent radical $U_0$, the maximal split
torus $T_0$, and the center $Z_0$. The restriction of scalars
gives corresponding subgroups $B/\Q, T/\Q, U/\Q$ and $Z/\Q$ of
$G$. We single out the element $w_0=\begin{pmatrix}
  0 & 1 \\
  -1 & 0
\end{pmatrix} \in G(\Q)$.

The positive simple root defines a homomorphism
$$\alpha_0:B_0/F \to \mathbf{G}_m/F$$
$$ \begin{pmatrix}
  t_{1} & * \\
  0 & t_{2}
\end{pmatrix}
\mapsto t_1/t_2$$
 and we denote by $\alpha$ the corresponding homomorphism
 $B/\Q \to {\rm Res}_{F/\Q} \mathbf{G}_m$. From \cite{HaGL2} we take the notation
 $|\alpha|$ for $| \, | \circ \alpha_{\A}:B(\A) \to \C^*$, where
 $| \, |: F^* \backslash \A_F^* \to \C^*$ is the idelic absolute value $x
 \mapsto |x|= \prod_v |x_v|_v$. Here we take the usual normalized
 absolute values for the local absolute values, in particular, $|x_{\infty}|_{\infty}=x_{\infty}
 \overline{x}_{\infty}$ at the complex place.

In $G_{\infty}={\rm GL}_2(\C)$ we choose the subgroup
$K_{\infty}=U(2) \cdot Z_0(\C)=U(2) \cdot \C^*$ containing the
maximal compact subgroup of unitary matrices. The symmetric space
$X=G_{\infty}/K_{\infty}$ can be identified with
three-dimensional hyperbolic space $\h=\Real_{>0} \times \C$.

The Lie algebra $\mathfrak{g}={\rm Lie}(G/\Q)$ is a $\Q$-vector
space and we define $\mathfrak{g}_{\infty}=\mathfrak{g}
\otimes_{\Q} \Real$. It carries a positive semidefinite
$K_{\infty}$-invariant form, the Killing form
$$\langle X,Y \rangle=\frac{1}{16} {\rm trace} ({\rm ad}X \cdot {\rm ad} Y),$$
and with respect to this form we have an orthogonal decomposition
$\mathfrak{g}_{\infty}=\mathfrak{k}_{\infty} \oplus \mathfrak{p}$,
where $\mathfrak{k}_{\infty}={\rm Lie}(K_{\infty})$ and
$$\mathfrak{p}=\Real H \oplus \Real E_1 \oplus \Real E_2:=\Real
\begin{pmatrix} 1&0\\0&-1\end{pmatrix} \oplus \Real
\begin{pmatrix} 0&1\\1&0\end{pmatrix} \oplus \Real
\begin{pmatrix} 0&i\\-i&0\end{pmatrix}.$$
Put $$S_{\pm}:= 1/2 \left( \pm \begin{pmatrix} 0&1\\1&0
\end{pmatrix} \otimes_{\Real} 1 -
\begin{pmatrix} 0&i\\-i&0 \end{pmatrix} \otimes_{\Real} i \right) \in
\mathfrak{p}_{\C}.$$

A maximal open compact subgroup of $G(\A_f)$ is given by $${\rm
GL}_2(\widehat{\mathcal{O}})=\left \{
\begin{pmatrix}
  a & b \\
  c & d
\end{pmatrix}
:a,b,c,d \in \widehat{\mathcal{O}}, ad-bd \in
\widehat{\mathcal{O}}^* \right \} .$$ We will deal with the
following congruence subgroups: For an ideal $\mathfrak{N}$ in
$\mathcal{O}$ and a finite place $v$ of $F$ put
$\mathfrak{N}_v=\mathfrak{N} \mathcal{O}_v$. We define

$$K^1(\mathfrak{N})=\left \{
\begin{pmatrix}
  a & b \\
  c & d
\end{pmatrix}
\in {\rm GL}_2(\widehat{\mathcal{O}}), a-1, c \equiv 0 \,
\mod{\mathfrak{N}}  \right \},$$

$$K^1(\mathfrak{N}_v)=\left \{
\begin{pmatrix}
  a & b \\
  c & d
\end{pmatrix}
\in {\rm GL}_2(\mathcal{O}_v), a-1, c \equiv 0 \,
\mod{\mathfrak{N}_v} \right \},$$ and $$U^1(\mathfrak{N}_v)=\{k
\in {\rm GL}_2(\mathcal{O}_v): {\rm det}(k) \equiv 1
\mod{\mathfrak{N}_v} \}.$$

For any compact open subgroup $K_f \subset G(\A_f)$ the adelic
symmetric  space $$S_{K_f}:=G(\Q)\backslash G(\A) /K_{\infty}
K_f$$  has several connected components. In fact, strong
approximation implies that the fibers of the determinant map
$$S_{K_f}
\twoheadrightarrow \pi_0(K_f):=\A_{F,f}^*/{\rm det}(K_f)F^*$$ are
connected. Any $\gamma \in G(\A_f)$ gives rise to an injection
$j_{\gamma}: G_{\infty} \to G(\A)$ with $j_{\gamma}(g_{\infty}):=
(g_{\infty},\gamma)$ and, after taking quotients, to a component
$$\Gamma_{\gamma}\backslash G_{\infty}/K_{\infty} \to S_{K_f},$$
where $\Gamma_{\gamma}:= G(\Q) \cap \gamma K_f \gamma^{-1}$. This
component is the fiber over ${\rm det}(\gamma)$. Choosing a system
of representatives for $\pi_0(K_f)$ we therefore have
$$S_{K_f} \cong \coprod_{[{\rm det}(\gamma)] \in \pi_0(K_f)}
\Gamma_{\gamma} \backslash \h.$$ We denote the Borel-Serre
compactifications of $S_{K_f}$ and $\Gamma_{\gamma} \backslash \h$
by $\overline S_{K_f}$ and $\Gamma_{\gamma} \backslash \oh$,
respectively. Following \cite{BS} we write $e(P)=\h/A_P\cong
U_P(\Real)$ for each rational Borel subgroup $P$ of $G$. Here
$U_P$ denotes its unipotent radical and $A_P$ the identity
component of $P(\Real)/U_P(\Real)$, and the action of $A_P$ on
$\h$ is the geodesic action. The boundary of $\Gamma_{\gamma}
\backslash \oh$ is the union of tori $\Gamma_{\gamma, P}
\backslash e(P)=:e'(P)$ with $\Gamma_{\gamma,P}=\Gamma_{\gamma}
\cap P(\Q)$ over a set of representatives for the
$\Gamma_{\gamma}$-conjugacy classes of Borel subgroups
(equivalently of $B(\Q) \backslash G(\Q)/\Gamma_{\gamma} \cong
\mathbf{P}^1(F)/\Gamma_{\gamma}$). We recall from \cite{HaGL2} \S
2.1 and \cite{Ha82} p. 110 that $\partial \overline S_{K_f}$ is
homotopy equivalent to \begin{equation} \label{BSbdry} \partial
\tilde S_{K_f}:=B(\Q) \backslash G(\A)/K_f K_{\infty} \cong
\coprod_{[{\rm det}(\gamma)] \in \pi_0(K_f)} \coprod_{[\eta] \in
\mathbf{P}^1(F)/\Gamma_{\gamma}} \Gamma_{\gamma, B^{\eta}}
\backslash \h,\end{equation} where $B^{\eta}(\Q)=\eta^{-1} B(\Q)
\eta $ for $\eta \in G(\Q)$ and the boundary component
$\Gamma_{\gamma, B^{\eta}} \backslash \h$ gets embedded in
$\partial \tilde S_{K_f}$ via $g_{\infty} \mapsto j_{\eta,
\gamma}(g_{\infty}):=\eta (g_{\infty}, \gamma)$.

\subsection{Hecke characters} \label{Hcharacters}
 A Hecke character of $F$ is a continuous group homomorphism
 $\lambda: F^* \backslash \A_F^* \to \C^*$. Such a character corresponds
 uniquely to a character on ideals prime to the conductor (see
 \cite{Hi93} \S 8.2), which we will also denote by $\lambda$.
  The archimedean part
 $\lambda_{\infty} : \C^* \to \C^*$ is of the form $z \mapsto
 \frac{z^a \overline{z}^b}{(z \overline{z})^t}$ for $t \in \C, a,b \in \Z$.
 We will say that $\lambda$ has \textit{infinity type}
 $\frac{z^a \overline{z}^b}{(z \overline{z})^t}$.
We define the (incomplete) $L$-series $L(s, \lambda)$ for ${\rm
Re}(s) \gg 0$ by the Euler product

$$L(s, \lambda):= \prod_{v \nmid \mathfrak{f}_{\lambda}} (1- \lambda(\mathfrak{P}_v)
{\rm Nm}(\mathfrak{P}_v)^{-s})^{-1},$$ where
$\mathfrak{f}_{\lambda}$ is the conductor of $\lambda$. This can
be continued to a meromorphic function on the whole complex plane
and satisfies a functional equation (see e.g., \cite{Hi93} \S 8.6
or \cite{dS} 37).

Define the character $\lambda^c$ by
$\lambda^c(x)=\lambda(\sigma(x))$. Since $\sigma$ just permutes
the Euler factors we have $L(s, \lambda)=L(s, \lambda^c)$. Also
let $\lambda^*(x):=\lambda(\sigma(x))^{-1} |x|$.

Recall from \cite{dS} p.91 and \cite{La} XIV Theorem 14 the
definition of the global root number $W(\lambda)$ appearing in the
functional equation. Note that $W(\lambda)=W(\tilde \lambda)$ for
$\tilde \lambda$ the associated unitary character
$\lambda/|\lambda|$. If $\lambda^*=\lambda$ then one shows using
the functional equation that $W(\lambda)= \pm 1$. For $\lambda$ of
infinity type $\frac{z^m}{(z \overline z)^{m/2}}$ with $m \in \Z$
we have
$$W(\lambda)=i^{-m} ({\rm Nm}(\mathfrak{f}_{\lambda}))^{-1/2} \prod_{v
\mid \mathfrak{f}_{\lambda}} \tau_v(\lambda) \prod_{v \nmid
\mathfrak{f}_{\lambda}} \lambda(\mathcal{D}_v^{-1}),$$ where the
Gauss sum $\tau_v$ is given by
$$\tau_v(\lambda_v)=\sum_{\epsilon \in
\mathcal{O}_v^*/(1+f_{\lambda,v})} (\lambda \mathbf{e}_F)
(\epsilon \pi^{-{\rm ord}_v(\mathfrak{f}_{\lambda}
\mathcal{D})}).$$ Here $\mathbf{e}_F$ is the standard additive
character of $F \backslash \A_F$ defined by $e_F=e_{\Q} \circ {\rm
Tr}_{F/\Q}$ in terms of the standard additive character $e_{\Q}$
of $\Q \backslash \A$ normalized by $e_{\Q}(x_{\infty})=e^{2 \pi i
x_{\infty}}$. Put $\tau(\lambda)=\prod_{v \mid
\mathfrak{f}_{\lambda}} \tau_v(\lambda)$.

We will use the following formula of Weil as stated in \cite{TA}
Proposition 2.4:

\begin{prop} \label{rootprod}
Suppose that $\lambda_1$ and $\lambda_2$ are unitary Hecke
characters of infinity types $(k_1, j_1)$ and $(k_2, j_2)$ with
relatively primes conductors $\mathfrak{f}_1$ and
$\mathfrak{f}_2$. Then

$$W(\lambda_1) W(\lambda_2) \lambda_1(\mathfrak{f}_2)
\lambda_2(\mathfrak{f}_1)=
  \begin{cases}
    W(\lambda_1 \lambda_2) & \text{if } (k_1-j_1)(k_2-j_2) \geq 0, \\
    (-1)^{\nu}W(\lambda_1 \lambda_2) & \text{if } (k_1-j_1)(k_2-j_2) < 0,
  \end{cases}$$
where $\nu={\rm min} \{|k_1-j_1|,|k_2-j_2| \}$. \hspace{\fill}
\qedsymbol
\end{prop}

For ease of reference we record the following:
\begin{lem} \label{charinteg}
For $\lambda: F^* \backslash \A_F^* \to \C^*$ with infinity type
$z^a \overline z^b$ with $a, b \in \Z$ we denote by
$\mathcal{O}_{\lambda}$ the ring of integers in the finite
extension of $F_{\mathfrak{p}}$ obtained by adjoining the values
of the finite part of $\lambda$. Then for any $x \in \A_{F,f}^*$
$${\rm ord}_{\mathfrak{p}}(\lambda(x))=-a \cdot {\rm ord}_{\mathfrak{p}}(x_{\mathfrak{p}})- b \cdot {\rm ord}_{\overline{\mathfrak{p}}}(x_{\overline{\mathfrak{p}}}).$$
\end{lem}
\begin{proof}
Let $v$ be any finite place of $F$. Since $\lambda$ has finite
order on $\mathcal{O}_v^*$ it suffices to prove the statement for
$\lambda(\pi_v)$ for any uniformizer $\pi_v$. If $h$ is the class
number of $F$, we have $\mathfrak{P}_v^h=(\alpha)$ for $\alpha \in
\mathcal{O}$ and $\alpha \in \mathcal{O}_w^*$ for $w \neq v$. Now
$$1=\lambda((\alpha, \alpha, \ldots))=\lambda_{\infty}(\alpha)
\lambda_v(\alpha) \prod_{w \neq v} \lambda_w(\alpha).$$ Since
$\prod_{w \neq v} \lambda_w(\alpha) \in \mathcal{O}_{\lambda}^*$
we deduce that $$h \cdot {\rm
ord}_{\mathfrak{p}}(\lambda(\pi_v))={\rm
ord}_{\mathfrak{p}}(\lambda_v(\alpha))=- {\rm
ord}_{\mathfrak{p}}(\lambda_{\infty}(\alpha)).$$
\end{proof}

Define $\Omega \in \C$ to be the complex period of a N\'{e}ron
differential $\omega$ of an elliptic curve $E$ defined over some
number field such that $E$ has complex multiplication by
$\mathcal{O}$, $E$ has good reduction at the place above $p$ and
$\overline \omega$ is a non-vanishing invariant differential on
the reduced curve $\overline E$.

Let $\lambda$ be a  Hecke character of infinity type $z^a
\overline z^b$ with $a,b \in \Z$. Precisely for $a>0$ and $b \leq
0$ or $a \leq 0$ and $b>0$ the $L$-value $L(0,\lambda)$ is
critical in the sense of Deligne. Damerell showed in this case
that $\pi^{{\rm max}(-a,-b)} \Omega^{-|a-b|}L(0,\lambda)$ is an
algebraic number in $\C$. We recall the following results (due to,
amongst others, Shimura, Coates-Wiles, Katz, Hida, Tilouine, de
Shalit, and Rubin) about the integrality of the special $L$-value
at $s=0$:
\begin{thm}  \label{Linteg}
Let $\lambda$ a Hecke character of infinity type $z^a \overline
z^b$ with conductor prime to $p$. Assume $a,b \in \Z$ and $a>0$
and ${b} \leq 0$. Put
$$L^{\mathrm{alg}}(0,\lambda):= \Omega^{b-a} \left(
\frac{2\pi}{\sqrt{d_F}}\right)^{-b} \Gamma(a)\cdot L(0,\lambda).$$

  \emph{(a)} If $p$ is split then
$$ (1-\lambda(\overline{\mathfrak{p}}))(1-\lambda^*(\overline{\mathfrak{p}})) \cdot
L^{\mathrm{alg}}(0,\lambda)$$ lies in the ring of integers of a
finite extension of $F_{\mathfrak{p}}$.

  \emph{(b)} If $p$ is inert and $a>0, b=0$
  then for any  ideal $\mathfrak{b}$ coprime to 6p and the conductor of $\lambda$
$$({\rm Nm}(\mathfrak{b})-\lambda^{-1}(\mathfrak{b})) \cdot L^{\mathrm{alg}}(0,\lambda)$$
lies in the ring of integers of a finite extension of
$F_{\mathfrak{p}}$.
\end{thm}

\begin{proof}[References]
If $p$ is split then the normalization in (a) is the one appearing
in the $p$-adic $L$-function constructed by Manin-Vishik, Katz,
and others. Together, \cite{K76} Chapters 4 and 8, \cite{K78}
Theorem 5.3.0, and \cite{HT} Theorem II prove that it is a
$p$-adic integer in $\widehat{\overline F}_{\mathfrak{p}}$. With
our fixed embedding $\overline F \hookrightarrow \overline
F_{\mathfrak{p}}$ this shows that the value lies in a finite
extension of $F_{\mathfrak{p}}$ and is $p$-integral. See also
\cite{Hi04a} Theorem 1.1 and \cite{dS} Theorem II.4.14 and II.6.7.

 Part (b) uses
the relation of elliptic units to special values of $L$-functions.
For the proof in the case when $\lambda$ is the power of a
Gr\"ossencharacter of a CM elliptic curve and $F$ has class number
one see, for example \cite{RuCime} \S7, in particular, Theorem
7.22. To extend to the general case use the arguments in \cite{dS}
Chapter II.
\end{proof}

\begin{rem}
\begin{enumerate}
\item If $p$ is split then Lemma \ref{charinteg} shows that for $a \geq 2$
the factor $(1- \lambda^*(\overline{\mathfrak{p}}))$ is a
$p$-unit.

  \item
 If $F$ has class number one, $p>a$, and $\lambda$ is the power of a
Gr\"ossencharacter of a CM elliptic curve then \cite{Dee} Lemma
3.4.5 proves that there always exists an ideal $\mathfrak{b}$ such
that ${\rm Nm}(\mathfrak{b})-\lambda^{-1}(\mathfrak{b})$ is prime
to $p$.

\item For completeness we want to mention that for inert primes $p$ additional divisibilities have been
obtained in \cite{K77}, \cite{K82}, \cite{Ru83}, \cite{Fuj}, and
\cite{Chel}.

\end{enumerate}
\end{rem}

\subsection{Modules and Sheaves} \label{mods}
The group ${\rm GL}_2(F)$ acts on the $F$-vector space $M^n:= {\rm
Sym}^n(F^2)$ of homogeneous polynomials of degree $n$ in two
variables $X$ and $Y$ with coefficients in $F$ by right
translation:
\begin{equation*} \label{moebius} \begin{pmatrix}
  a & b \\
  c & d
\end{pmatrix}
 . X^iY^{n-i} = (aX+cY)^i (bX+dY)^{n-i}.\end{equation*} Applying first the field
automorphism $\sigma$ to the entries $a,b,c$ and $d$, we get
another representation $\overline{M}^n$. We also have
one-dimensional representations $F[k,\ell]$ for $(k,\ell) \in
\Z^2$, on which $g\in G$ acts by multiplication by ${\rm det}^k(g)
\cdot \sigma({\rm det}(g))^\ell$. We obtain the representations
$M(m,n,k,\ell):=M^m \otimes_F \overline{M}^n \otimes_F F[k,\ell]$.
Let $M(m,n,k,\ell)^{\vee}:={\rm Hom}_F(M(m,n,k,\ell), F).$ There
is an isomorphism of ${\rm GL}_2(F)$-modules
$$M(m,n,k,\ell)^{\vee} \cong M(m,n,-m-k,-n-\ell)$$
induced by the pairing
$$\langle \, , \, \rangle : M(m,n,k,\ell) \times M(m,n,-m-k,-n-\ell) \to F,$$
$$X^jY^{m-j}\overline X^k \overline Y^{n-k} \times X^{\mu}Y^{m-\mu}\overline X^{\nu} \overline
Y^{n-\nu} \mapsto (-1)^{j+k}
  \begin{pmatrix}
    m \\
    j
  \end{pmatrix}^{-1}
  \begin{pmatrix}
    n \\
    k
  \end{pmatrix}^{-1}
\delta_{j, m-\mu} \delta_{k,n-\nu}.$$ This is the coordinatized
version of the pairing induced by the determinant pairing on $F^2$
(cf. \cite{Hi93} p. 169).

For an $\mathcal{O}$-module $N$ we denote $N \otimes_{\mathcal{O}}
A$ by $N_A$ for any $\mathcal{O}$-algebra $A$. Denote by
$M(m,n,k,\ell)_{\mathcal{O}}$ the polynomials with
$\mathcal{O}$-coefficients. Note that
$M(m,n,k,\ell)_{\mathcal{O}}^{\vee}:={\rm
Hom}_{\mathcal{O}}(M(m,n,k,\ell), \mathcal{O})$ corresponds under
the duality above to $$\left\{\sum_{\mu,\nu} a_{\mu,\nu}
\begin{pmatrix}
    m \\
    \mu
  \end{pmatrix}
  \begin{pmatrix}
    n \\
    \nu
  \end{pmatrix}
  X^{\mu}Y^{m-\mu} \overline X^{\nu} \overline Y^{n-\nu}|
  a_{\mu,\nu} \in \mathcal{O}
\right \} \subset M(m,n,k,\ell).$$

We now define local coefficient systems on the symmetric spaces.
For $\Gamma \subset G(\Q)$ an arithmetic subgroup and $N$ an
$\mathcal{O}[\Gamma]$-module we define a sheaf of
$\mathcal{O}$-modules on $\Gamma \backslash \h$ by
\begin{align}\widetilde N(U):=
 \{ &f:\pi_{\Gamma}^{-1}(U) \to N \text{ locally constant }:
 \nonumber
 \\ \nonumber & f(\beta
x)=\beta.f(x) \forall x \in \pi_{\Gamma}^{-1}(U) \text{ and }
\beta \in \Gamma \},\nonumber \end{align} where $\pi_{\Gamma}: \h
\to \Gamma \backslash \h$ is the canonical projection.

Let $K_f \subset G(\A_f)$ be a compact open subgroup and $M$ an
$F[G(\Q)]$-module. Assume that there exists an
$\mathcal{O}$-lattice $M_{\mathcal{O}}$ in $M$ such that
$M_{\hat{\mathcal{O}}}=M_{\mathcal{O}} \otimes \hat{\mathcal{O}}$
is stable under $K_f$. (For $M=M(m,n,k,\ell)$ and $K_f \subset
{\rm GL}_2(\hat{\mathcal{O}})$ one can take
$M_{\mathcal{O}}=M(m,n,k,\ell)_{\mathcal{O}}$.) For each open
subset $U \subset S_{K_f}$ we let
$$ \widetilde M_{\mathcal{O}}(U):=\left\{ f: \pi^{-1}(U) \to M \text{ locally constant }
 \: \left| \begin{aligned}  &f(\beta g)= \beta.f(g), f(g) \in g_f M_{\hat{\mathcal{O}}}
 \\ &  \forall g \in \pi^{-1}(U)\text{ and } \beta \in G(\Q) \end{aligned} \right.
\right\},$$ where $\pi:G(\A)/K_{\infty}K_f \to S_{K_f}$ is the
 projection. This defines a sheaf of $\mathcal{O}$-modules on $S_{K_f}$
 (cf. \cite{U98} \S 1.4, \cite{Ko} \S 1.5, and \cite{F} \S 1.2).
For any $\mathcal{O}$-algebra $R$ we define $\widetilde M_R$ as
$\widetilde M_{\mathcal{O}} \otimes \underline{R}$, where
$\underline{R}$ is the constant sheaf associated to $R$.

For $\gamma \in G(\A_f)$ let $M_{\gamma} := M \cap
\gamma.M_{\hat{\mathcal{O}}}$. Then $M_{\gamma}$ is a locally
free, finitely generated $\mathcal{O}$-module with an action by
$\Gamma_{\gamma}=G(\Q) \cap \gamma K_f \gamma^{-1}$. The two
constructions of $\widetilde M_{\mathcal{O}}$ and $\widetilde
M_{\gamma}$ are compatible with $j_{\gamma}$; one checks that
$j_{\gamma}^*(\widetilde M_{\mathcal{O}}) \cong \widetilde
M_{\gamma}$.

\subsection{Cohomology} \label{Scohom}
For a sheaf $\mathcal{F}$ on a topological space $X$, we denote by
$H^i(X,\mathcal{F})$ (resp. $H^i_c(X,\mathcal{F})$) the $i$-th
cohomology group of $\mathcal{F}$ (resp. with compact support),
and the interior cohomology, i.e., the image of
$H^i_c(X,\mathcal{F})$ in $H^i(X,\mathcal{F})$, by
$H^i_!(X,\mathcal{F})$.

Let $M$ be an $F[G(\Q)]$-module with $M_{\mathcal{O}} \subset M$
an $\mathcal{O}$-lattice as above and $R$ an
$\mathcal{O}$-algebra. Since $S_{K_f} \overset {i}
{\hookrightarrow} \overline S_{K_f}$ is a homotopy equivalence, we
have a canonical isomorphism
$$H^i(S_{K_f}, \widetilde{M}_R) \cong H^i(\overline S_{K_f}, i_*
\widetilde{M}_R)$$ and in what follows we will replace $i_*
\widetilde{M}_R$ by $\widetilde{M}_R$ and also write
$\widetilde{M}_R$ for the sheaf $j^*i_* \widetilde{M}_R$ on
$\partial \overline S_{K_f}$, for $j:\partial \overline S_{K_f}
\hookrightarrow \overline S_{K_f}$.

The decomposition of the adelic symmetric space into connected
components gives rise to canonical isomorphisms (see \cite{Ko} \S
1.6 and \cite{F} \S 1.2)
$$H^i(S_{K_f}, \widetilde M_{R}) \cong \bigoplus_{[{\rm det}(\gamma)] \in \pi_0(K_f)}
H^i(\Gamma_{\gamma} \backslash \h, \widetilde{M_{\gamma}}\otimes
\underline R)$$ and
$$H^i(\partial \tilde S_{K_f}, \widetilde M_{R}) \cong \bigoplus_{[{\rm det}(\gamma)] \in
\pi_0(K_f)}\bigoplus_{[\eta] \in \mathbf{P}^1(F)/\Gamma_{\gamma}}
H^i(\Gamma_{\gamma, B^{\eta}} \backslash \h,
\widetilde{M_{\gamma}}\otimes \underline R).$$ The above
cohomology groups and isomorphisms are all functorial in $R$.

For an arithmetic subgroup $\Gamma \subset G(\Q)$ and
 an $\mathcal{O}[\Gamma]$-module $N$ we can
in many cases relate the sheaf cohomology $H^i(\Gamma \backslash
\h, \widetilde N_R)$ to group cohomology $H^i(\Gamma, N_R)$ (for
the proof see, e.g., \cite{HaCAG}):
\begin{prop} \label{shgp}
For $\mathcal{O}$-algebras $R$ in which the orders of all finite
subgroups of $\Gamma$ are invertible there is a natural
$R$-functorial isomorphism
$$H^i(\Gamma \backslash \h, \widetilde N_R) \cong H^i(\Gamma,
N_R).$$\hfill \qedsymbol
\end{prop}
The lemma in \cite{F} \S 1.1 shows that for any
$\mathcal{O}$-algebra $R$, $R\otimes_{\mathcal{O}}
\mathcal{O}[\frac{1}{6}]$  satisfies the conditions of the
proposition for any arithmetic subgroup $\Gamma \subset G(\Q)$.

For complex coefficient systems we have analytic tools available.
For a $C^{\infty}$-manifold $X$ (like $S_{K_f}$, $\partial \tilde
S_{K_f}$, or $\Gamma \backslash \h$) denote by $\Omega^i(X)$ the
space of $\C$-valued $C^{\infty}$-differential $i$-forms and by
$\Omega^i(X,M_{\C})=\Omega^i(X) \otimes_{\C} M_{\C}$ the space of
$M_{\C}$-valued smooth $i$-forms. By the de Rham Theorem (cf.
\cite{HaAG} IV.9.1, or \cite{Hi93} Appendix Theorem 2) we have
$$H^i(\Gamma\backslash \h, \widetilde {M_{\C}}) \cong
H^i(\Omega^{\bullet}(\h, M_{\C})^{\Gamma}).$$ Furthermore, the de
Rham cohomology groups are canonically isomorphic to relative Lie
algebra cohomology groups. For the definition of the latter we
refer to \cite{BW} Chapter 1. The tangent space of $\h$ at the
point $K_{\infty} \in G_{\infty}/K_{\infty}$ can be canonically
identified with $\mathfrak{g}_{\infty}/\mathfrak{k}_{\infty}$. For
$g \in G_{\infty}$ let ${\rm L}_g:\h \to \h$ be the
left-translation by $g$ and ${\rm D}_{{\rm L}_g}$ the differential
of this map. Assume that the $G(\Q)$-action on $M_{\C}$ extends to
a representation of $G_{\infty}$. Let $\omega_{M_{\C}}:Z(\Real)
\to \C^*$ be the character describing the action on $M_{\C}$ and
write $C^{\infty}(\Gamma \backslash {\rm
GL}_2(\C))(\omega_{M_{\C}}^{-1})$ for those functions in
$C^{\infty}(\Gamma \backslash {\rm GL}_2(\C))$ on which
translation by elements in $Z(\Real)$ acts via
$\omega_{M_{\C}}^{-1}$.

We can then identify the $\C$-vector spaces
$$\Omega^i(\h, M_{\C})^{\Gamma} \cong {\rm
Hom}_{K_{\infty}}(\Lambda^i(\mathfrak{g}_{\infty}/\mathfrak{k}_{\infty}),
C^{\infty}(\Gamma \backslash {\rm GL}_2(\C))(\omega_{M_{\C}}^{-1})
\otimes M_{\C}),$$ by mapping an $M_{\C}$-valued differential form
$\tilde \omega$ to the $(\mathfrak{g}, K_{\infty})$-cocycle
$\omega$ given by $\omega(g)(\theta_1 \wedge \ldots \wedge
\theta_i):= g^{-1}. \tilde \omega(g K_{\infty})({\rm D}_{{\rm
L}_g}(\theta_1), \ldots, {\rm D}_{{\rm L}_g}(\theta_i))$. The
differentials of the complexes corresponds and we get (cf.
\cite{BW} VII Corollary 2.7)
$$H^i(\Gamma\backslash \h, \widetilde {M_{\C}}) \cong
H^i(\mathfrak{g}_{\infty},
    K_{\infty}, C^{\infty}(\Gamma \backslash {\rm GL}_2(\C))
    (\omega^{-1}_{{M_{\C}}})
    \otimes {M_{\C}}).$$
Similarly, one obtains
$$H^i(S_{K_f}, \widetilde{{M_{\C}}}) \cong H^i(\mathfrak{g}_{\infty},
    K_{\infty}, C^{\infty}(G(\Q) \backslash G(\A)/K_f)(\omega^{-1}_{{M_{\C}}})
    \otimes {M_{\C}}) $$ and
$$H^i(\partial \tilde S_{K_f}, \widetilde{{M_{\C}}}) \cong H^i(\mathfrak{g}_{\infty},
    K_{\infty}, C^{\infty}(B(\Q) \backslash G(\A)/K_f)(\omega^{-1}_{{M_{\C}}})
    \otimes {M_{\C}}).$$
For any cocycle $\omega$ we will denote by $[\omega]$ the
corresponding cohomology class.

For $\Gamma \backslash \h$ and $N$ an $\mathcal{O}[\Gamma]$-module
the natural isomorphisms of the de Rham Theorem and Proposition
\ref{shgp} compose to give an isomorphism between de Rham
cohomology and group cohomology. We state this isomorphism
explicitly on the level of cocycles for degree 1 (for a proof see
\cite{Be} Proposition 2.5 or, more generally, \cite{BC} Proof of
Lemma 3.3.5.1):
\begin{prop} \label{dRgroup}
The natural isomorphism
$$H^1(\Omega^{\bullet} (\h,
{N_{\C}})^{\Gamma}) \cong H^1(\Gamma \backslash \h, \widetilde
{N_{\C}}) \cong H^1(\Gamma, {N_{\C}})$$  is induced by any of the
following maps on closed 1-forms:  For a choice of basepoint $x_0
\in \h$ assign to a closed 1-form $\tilde \omega$ with values in
${N_{\C}}$ the (inhomogeneous) 1-cocycle
$$\mathcal{G}_{x_0}(\tilde \omega): \alpha \mapsto \int_{x_0}^{\alpha.x_0}
\tilde \omega.$$ \hfill \qedsymbol
\end{prop}

For each $g \in G(\A_f)$ with $g \widetilde M_{\mathcal{O}}
\subset \widetilde M_{\mathcal{O}}$ we have the Hecke algebra
action of the double coset $[K_f g K_f]$ on the cohomology groups
$H^i(S_{K_f}, \widetilde M_R)$, $H^i(\partial \overline S_{K_f},
\widetilde M_R)$, and $H^i_!(S_{K_f}, \widetilde M_R)$ for any
$\mathcal{O}$-algebra $R$ (for its definition see \cite{U98}
\S1.4.4).

For $R=\C$ this can be described on the level of relative Lie
algebra cohomology: If $V$ is any $G(\A_f)$-module then $[K_f g
K_f]$ acts on a $K_f$-invariant vector $v \in V$ by $$[K_f g
K_f].v= \sum_{\gamma \in K_f g K_f/K_f} \gamma.v .$$ Taking this
action on $C^{\infty}(G(\Q) \backslash G(\A)/K_f)$,
$C^{\infty}(B(\Q) \backslash G(\A)/K_f)$, and \linebreak
$C^{\infty}_c(G(\Q) \backslash G(\A)/K_f)$, respectively, induces
via relative Lie algebra cocycles the Hecke action on the
cohomology groups. For $x \in \mathcal{O} \otimes \hat \Z$ we
single out the operators
$$T_x=[K_f \begin{pmatrix} x &0\\0&1 \end{pmatrix} K_f].$$

\subsection{Relative (co-)homology} \label{relativecoh}

We refer to \cite{B67} II \S 12 and V \S 5 for the definitions of
relative sheaf cohomology and relative Borel-Moore homology, but
want to recall the following facts:

Let $\Gamma \subset G(\Q)$ be an arithmetic subgroup , $K_f
\subset G(\A_f)$ a compact open subgroup, and $M$ an
$F[G(\Q)]$-module with $M_{\mathcal{O}} \subset M$ an
$\mathcal{O}$-lattice as above. Then for $X=\Gamma \backslash \oh$
or $\overline S_{K_f}$ and $Y \subset X$ a closed subspace we have
the long exact sequence (functorial in the $\mathcal{O}$-algebra
$R$)
$$ \ldots \to H^i(X,Y,\widetilde M_R) \to H^i(X,\widetilde M_R)
\to H^i(Y,\widetilde M_R) \to H^{i+1}(X,Y,\widetilde M_R) \to
\ldots$$ Note that for $Y=\partial \overline S_{K_f}$ we have
$H^i(X,Y,\widetilde M_R) \cong H^i_c(X, \widetilde M_R)$.

For $Y \subset \partial \overline S_{K_f}$ we get the following
analytic description: We say that a function $f \in
 C^{\infty}(\Gamma \backslash {\rm GL}_2(\C))$ has
\textit{moderate growth} if there exists an integer $N>0$ and a
constant $c>0$ such that
$$|f(g)| \leq c \|g\|^N, \, \text{ for all } g \in {\rm
GL}_2(\C).$$ For any Borel subgroup $P$ of ${\rm GL}_2$ defined
over $F$ we say that a function $f \in C^{\infty}(\Gamma
\backslash {\rm GL}_2(\C))$ is \textit{fast decreasing at $P$} if
$f$ has moderate growth and if there exists an integer $N>0$ such
that for every Siegel set $S \subset {\rm GL}_2(\C)$ relative to
$P$, every compact set $\omega \subset P(\Real)$, and every $r \in
\Real$ there exists a constant $c(S,\omega, r) >0$ satisfying
$$|f(a\omega g k) | \leq c(S,\omega, r) \|a\|^N \|g\|^r,$$
for all $a \in Z, w \in \omega, g \in S \cap {\rm SL}_2(\C), k \in
K_{\infty}$.

Let $C(\Gamma)$ be the set of $\Gamma$-conjugacy classes of Borel
subgroups of ${\rm GL}_2$ defined over $F$. For $\mathfrak{c}
\subset C(\Gamma)$ denote by $C^{\infty}_{\mathfrak{c}}(\Gamma
\backslash {\rm GL}_2(\C))$ the space of functions $f$ which,
together with all their derivatives $Df , D \subset
U(\mathfrak{g})$, are of moderate growth and fast decreasing at
every $P$ such that $[P] \in \mathfrak{c}$. We then get the
following extension of the de Rham theory recalled in Section
\ref{Scohom}:
\begin{prop} \label{Borel} $$H^1( \Gamma \backslash \oh, \bigcup_{[P] \in \mathfrak{c}} e'(P),\widetilde M_{\C}) \cong
H^1(\mathfrak{g}_{\infty},K_{\infty},C^{\infty}_{\mathfrak{c}}(\Gamma
\backslash {\rm GL}_2(\C)) \otimes M_{\C}).$$ \end{prop}

\begin{proof}[Sketch of proof]
We need to show that the inclusion of the complex of differential
forms with compactly supported coefficients on $\Gamma \backslash
\oh - \left( \bigcup_{[P] \in \mathfrak{c}} e'(P) \right)$ into
the forms with fast decreasing coefficients induces an isomorphism
in cohomology. The corresponding statement for cohomology with
compact support (and for general locally symmetric spaces) is
proven by Borel in \cite{Bo} Theorem 5.2. Since the proof uses
sheaf theory and considers stalks in the boundary it extends to
our case of relative cohomology with respect to subsets of
$C(\Gamma)$.
\end{proof}

From the comparison theorem with relative singular cohomology (see
\cite{B67} X \S14) we obtain the evaluation pairing
$$H^i(X,Y,\widetilde M_R^{\vee})_{\rm free} \times H_i(X,Y,\widetilde
M_R)_{\rm free} \to R,$$ which is perfect for any
$\mathcal{O}[\frac{1}{6}]$-algebra $R$ (see \cite{F}, Satz 3, and
\cite{G67} \S23). For $R=\C$ this pairing can be calculated by
$([\omega], [\sigma]) \mapsto \int_{\sigma} \omega$ for $\omega$ a
relative Lie algebra cocycle and $\sigma$ a differentiable
singular cycle. Note that for $R \subset \C$ a class $[\omega] \in
H^i(X,Y,\widetilde M_{\C}^{\vee})$ lies in $H^i(X,Y,\widetilde
M_R^{\vee})_{\rm free} \cong {\rm im}(H^i(X,Y,\widetilde
M_R^{\vee}) \to H^i(X,Y,\widetilde M_{\C}^{\vee}))$ if and only if
all the pairings of $[\omega]$ with homology classes in
$H_i(X,Y,\widetilde M_R)_{\rm free}$ have values in $R$.

\section{Eisenstein cohomology}
In this section we recall Harder's construction of Eisenstein
cohomology, construct explicit classes, and calculate their Hecke
eigenvalues and restrictions to the boundary. We also investigate
the integrality of the latter by translating to group cohomology.
\subsection{Eisenstein cocycles} \label{s3.1} Let $m,n \in \N_{\geq 0}$, $k,\ell \in \Z$, and $M:=M(m,n, k,\ell)$. We
want to construct certain cohomology classes in $H^1(S_{K_f},
\tilde M_{\C})$ with nontrivial restriction to the boundary
following the work of Harder in \cite{Ha79}, \cite{Ha82}, and
\cite{HaGL2}.

Let $\phi_1, \phi_2:F^* \backslash \A_F^* \to \C^*$ be two Hecke
characters with conductors $\mathfrak{M}_1$ and $\mathfrak{M}_2$
and of infinity type either
   $$\phi_{1,\infty}(z)=z^{1-k} \overline z^{-n-\ell} \; \text{ and } \phi_{2,\infty}(z)=z^{-m-k-1} \overline z^{-\ell} \text{ (Case A)}$$
or
   $$\phi_{1, \infty}(z)=z^{-m-k} \overline z^{1-\ell} \; \text{ and } \phi_{2, \infty}(z)=z^{-k} \overline z^{-n-\ell-1}\text{ (Case B)}.$$
They determine a character $\phi=(\phi_1, \phi_2)$ on $T(\Q)
\backslash T(\A)$. We put $\chi:=\phi_1/\phi_2$ and denote its
conductor by $\mathfrak{M}$.
\begin{rem} These are the infinity types of Hecke characters contributing
to the cohomology of the boundary as calculated in \cite{HaGL2} \S
2.9 and 3.5. The two cases get swapped by the action of the Weyl
group, which is defined by $w_0 . (\phi_1,\phi_2)= (\phi_2  |
\cdot |, \phi_1 | \cdot |^{-1})$, and we are in the so-called
``balanced case" (cf. \cite{HaGL2} \S 2.9)
\end{rem}

For a continuous character $\eta:T(\Q) \backslash T(\A) \to \C^*$
we define the induced module
\begin{equation*} V_{\eta,\C}=
\left\{ \Psi:G(\A) \rightarrow \C
 \: \left| \begin{aligned}  &\Psi(bg)=\eta(b)\Psi(g), \:\forall b \in B(\A),
 \\ &\Psi(gk)=\Psi(g)   \: \forall k\in K_f \subset G(\A) \text{ compact open} \\
 &\Psi \, \text{is $K^{\infty}$-finite on the right }
\end{aligned} \right. \right\}.\end{equation*}
We use here the following convention: for any $\Q$-algebra $R$ we
consider characters $\eta$ of $T(R)$ as characters of $B(R)$ by
defining $\eta(b):= \eta(t)$ if $b=tu$ for $t \in T(R)$ and $u \in
U(R)$. Note that the definition for $V_{\eta}$ follows the one
used in Harder's work and is not the usual unitary induction. The
induced representation $V_{\eta,\C}$ decomposes into a product
$\bigotimes_v V_{\eta_v,\C}$, where $$V_{\eta_v, \C}=\{
\Psi_v:G_0(F_v) \to \C | \Psi_v(b_v g_v)= \eta_v(b_v) \Psi(g_v) \,
\forall b \in B_0(F_v) \}$$ denotes the local induced
representations. By $V_{\eta_f, \C }=\otimes_{v \nmid \infty}
V_{\eta_v, \C}$ we denote the finite part of $V_{\eta,\C}$. We
will be interested in the operation of the Galois group ${\rm
Gal}(\overline \Q/\Q)$ on $V_{\eta_f, \C}$ and write $V_{\eta_f}$
(resp. $V_{\eta_v}$  for $v\nmid \infty$) for the $\overline
\Q$-subspace of $V_{\eta_f, \C}$ (resp. $V_{\eta_v,\C}$)
consisting of $\overline \Q$-valued functions. Every $\sigma \in
{\rm Gal}(\overline \Q/\Q)$ defines a $\sigma$-linear isomorphism
$$ \sigma: V_{\eta_v} \to V_{\eta_v^\sigma}$$ $$\Psi \mapsto
\Psi^{\sigma},$$ where for each $g \in G_0(F_v)$ we define
$\Psi^{\sigma}(g):=\Psi(g)^{\sigma}$. By \cite{Wald}, ch. I.2 (see
also  \cite{Mahn} p. 94) we have $V_{\eta_v, \C}= V_{\eta_v}
\otimes \C$ which implies $V_{\eta_f,\C}=V_{\eta_f} \otimes \C$.
The Galois action on $V_{\eta_f}$ is defined similarly to that on
$V_{\eta_v}$.

Given $\phi$ of infinity type (A) or (B), $\Psi \in
V_{\phi_f}^{K_f}$, and an appropriate open compact subgroup $K_f
\subset G(\A_f)$ we will first define a boundary cohomology class
$$[\omega_0(\phi, \Psi)] \in H^1(\partial \overline S_{K_f},
\widetilde M_{\C})$$ and then an Eisenstein cohomology class
$$[{\rm Eis}(\phi, \Psi)] \in H^1(S_{K_f}, \widetilde M_{\C}).$$

Harder describes the cohomology of the boundary as a
$G(\A_f)$-module in Theorem 1 of \cite{HaGL2}:

\begin{equation}
H^1(\partial \overline S_{K_f}, \widetilde M_{\overline F}) \cong
\bigoplus_{
   \begin{array}{c}
    \phi:T(\Q) \backslash T(\A) \to \C^*\\
    \text{ of infinity type (A)}
   \end{array}}
 \left ( V_{\phi_f}^{K_f}
\oplus V_{w_0.\phi_f}^{K_f} \right).
\end{equation} By the proof of Theorem 2 of \cite{HaGL2} (see also
Proposition 2.12 of \cite{Be}) relative Lie algebra cocycles
giving rise to non-trivial cohomology classes in
$$H^1(\partial \overline S_{K_f}, \widetilde M_{\C}) \cong
H^1(\mathfrak{g}_{\infty},
 K_{\infty}, C^{\infty}(B(\Q) \backslash
 G(\A)/K_f)(\omega_{M_{\C}}^{-1}) \otimes M_{\C})$$ can be described by certain
elements in $${\rm
Hom}_{K_{\infty}}(\mathfrak{g}_{\infty}/\mathfrak{k}_{\infty},
V_{\phi,\C}^{K_f} \otimes M_{\C})\cong
 (\check{\mathfrak{p}}_{\C} \otimes_{\C} M_{\C} \otimes_{\C}
 V_{\phi,\C}^{K_f})^{K_{\infty}} $$ using the map in relative Lie algebra cohomology induced by the embedding
 $$V_{\phi,\C}^{K_f} \hookrightarrow C^{\infty}(B(\Q) \backslash
 G(\A)/K_f)(\omega_{M_{\C}}^{-1}).$$ Recall $|\alpha|:B(\A)\to \C^*$ from Section \ref{s2.2}. Following \cite{HaGL2} p. 80 and \cite{Ko} p. 101 we
 define $$\omega_z(\cdot, \phi, \Psi):G(\A) \to \check{\mathfrak{p}}_{\C} \otimes_{\C}
 M_{\C}$$ for $z\in \C$ and $\Psi \in V_{\phi_f |\alpha|_f^{z/2},\C}^{K_f}$  as

\begin{eqnarray} \label{omega0} &&\omega_z(g, \phi, \Psi):= \omega(b_{\infty} k_{\infty} \cdot
g_f, \phi |\alpha|^{z/2}, \Psi)=\\
&=& (\phi_{\infty} \cdot |\alpha|_{\infty}^{z/2})(b_{\infty})
\cdot \Psi(g_f)
  \begin{cases}
    k_{\infty}^{-1} . \,(\check{S}_+ \otimes Y^m \overline{X}^n ) & \text{ Case (A), } \\
    k_{\infty}^{-1} . \,((-\check{S}_-) \otimes X^m \overline{Y}^n ) & \text{ Case (B). }
  \end{cases}
  \nonumber \end{eqnarray}
Here $K_{\infty}$ acts on $\mathfrak{p}_{\C}$ by the adjoint
action. By \cite{Ha79} Lemma 1.5.2 $\omega_0$ is a relative Lie
algebra 1-cocycle. We write $[\omega_0(\phi, \Psi)]$ both for the
corresponding cohomology class in $H^1(\mathfrak{g}_{\infty},
\mathfrak{k}_{\infty}, V_{\phi} \otimes M_{\C})$ as well as its
non-trivial image in $H^1(\partial \overline S_{K_f},
\widetilde{M_{\C}})$.

We now have for ${\rm Re}(z) \gg 0$ an operator
$${\rm Eis}: V_{\phi_f |\alpha|_f^{z/2},\C}^{K_f} \to \mathcal{A}(G(\Q) \backslash G(\A)/K_f)$$
given by the formula $$\Psi \mapsto {\rm Eis}(\Psi)(g)=
\sum_{\gamma \in B(\Q) \backslash G(\Q)} \Psi(\gamma g).$$  This
can be meromorphically continued to all $z\in \C$. Via the map on
cocycles the operator induces a map in cohomology. Define ${\rm
Eis}(\phi |\alpha|^{z/2},\Psi):={\rm Eis} (\omega_z(\phi, \Psi))$
for $\Psi \in V_{\phi_f |\alpha|_f^{z/2},\C}^{K_f}$.

Harder shows in \cite{HaGL2} Theorem 2 that for $z=0$ we get a
holomorphic closed form. For $g \in G(\A)$ and $A \in
\mathfrak{g}_{\infty}/\mathfrak{k}_{\infty}$ we use the notation
${\rm Eis}(g,\phi, \Psi)(A)$ for the Lie algebra 1-cocycle in
${\rm
Hom}_{K_{\infty}}(\mathfrak{g}_{\infty}/\mathfrak{k}_{\infty},
C^{\infty}(G(\Q) \backslash G(\A)/K_f)(\omega^{-1}_{M_{\C}})
\otimes M_{\C})$. The corresponding de Rham 1-form in
$\Omega^1(G(\A)/K_f K_{\infty} \otimes M_{\C})^{G(\Q)}$ is denoted
by ${\rm Eis}(x, \phi, \Psi)(\theta_x)$ for $x \in S_{K_f}$ and
$\theta_x \in T_x S_{K_f}$. We write  $[{\rm Eis}(\phi, \Psi)]$
for the cohomology class in $H^1(S_{K_f}, \widetilde{M_{\C}})$. We
will later drop $\phi$ in the argument if it is clear from the
context.

\subsection{Special vectors} \label{s3.2} We now want to single out some special vectors
$$\Psi=\otimes_{v \nmid \infty} \Psi_v \in V_{\phi_f
|\alpha|_f^{z/2}}.$$ At finite places we define the following
functions:

\begin{enumerate}[(a)]
  \item For any finite place $v$ we define
  $\Psi^{\rm new}_v$ to be the newvector spanning
  $V_{\phi_v |\alpha|_v^{z/2}}^{K^1(\mathfrak{P}_v^s)}$,
  where $\mathfrak{P}_v^s \parallel \mathfrak{M}_1 \mathfrak{M}_2$.
  By \cite{Cas} \S1 $\mathfrak{P}_v^s$ is the conductor of $V_{\phi_v
  |\alpha|_v^{z/2}}$ and we normalize $\Psi^{\rm new}_v$ by:
\begin{equation*}\Psi_v^{\rm new}(g)=
  \begin{cases}
    \phi_{1,v}(a) \phi_{2,v}(d) \left | \frac{a}{d}\right|_v^{z/2} & \text{if } g=\begin{pmatrix} a & b\\ 0 & d
    \end{pmatrix}\begin{pmatrix} 1 & 0 \\ \pi_v^r & 1
    \end{pmatrix} k, k \in K^1(\mathfrak{P}_v^s) \\
    0 & \text{otherwise,}
  \end{cases}
\end{equation*}
where $\mathfrak{P}_v^r \parallel \mathfrak{M}_1$.

  \item For $v \nmid \mathfrak{M}$ we also have the spherical
  vector $\Psi^0_v \in V_{\phi_v
  |\alpha|_v^{z/2}}^{U^1(\mathfrak{M}_{1,v})}$ defined by
\begin{equation*} \Psi^0_v(g)= \phi_{1,v}(a) \phi_{2,v}(d) \left | \frac{a}{d}\right|_v^{z/2}
 \phi_{1,v}({\rm det}(k))
\text{ for } g=\begin{pmatrix} a & b\\ 0 & d
    \end{pmatrix} k, k\in {\rm GL}_2(\mathcal{O}_v).\end{equation*}
\end{enumerate}
Note that $\Psi_v^0=\Psi_v^{\rm new}$ for $v \nmid
\mathfrak{M}_1\mathfrak{M}_2$.

Denote by $S$ the finite set of places where both $\phi_i$ are
ramified, but $\chi=\phi_1/\phi_2$ is unramified. We put
$$\Psi_{\phi_f |\alpha|_f^{z/2}}^{\rm new}:=\prod_{v \nmid\infty} \Psi^{\rm new}_v
\in V_{\phi_f |\alpha|_f^{z/2}}^{K_f^{\rm new}}$$
$$\Psi_{\phi_f |\alpha|_f^{z/2}}^0:=\prod_{v \notin S,v \nmid \infty}
\Psi_v^{\rm new} \prod_{v \in S} \Psi_v^0 \in V_{\phi_f
|\alpha|_f^{z/2}}^{K_f^S},$$ where $$K^{\rm
new}_f:=K^1(\mathfrak{M}_1 \mathfrak{M}_2)$$ and
$$K^S_f:=\prod_{v \in S}
U^1(\mathfrak{M}_{1,v}) \prod_{v \notin S} K^1((\mathfrak{M}_1
\mathfrak{M}_2)_v).$$

The following lemmata from \cite{Be} tell us how to translate
between $\Psi^0_{v}$ and $\Psi^{\rm new}_{v}$: Let
$\eta=(\eta_1,\eta_2): T(\Q)\backslash T(\A) \to \C^*$ be a
continuous character, e.g. $\eta=\phi |\alpha|^{z/2}$.

\begin{lem}[(\cite{Be} Lemma 4.3)] \label{l3.2}
Let $v$ be a place where both $\eta_i$ are unramified, and
$\mu:F_v^* \to \C^*$ a continuous character. If $\Psi^{\rm new}
_{\eta_v}$ is the newvector in $V_{\eta_v}$ and $\Psi^0_{\eta_v
\mu}$ is the spherical vector in $V_{\eta_v \mu}$ then $\Psi^{\rm
new} _{\eta_v}(g) \mu({\rm det}(g))=\Psi^0_{\eta_v \mu}(g)$ for
all $g \in F_v^*$. \hfill \qedsymbol
\end{lem}

\begin{lem} [(\cite{Be} Lemma 4.4)] \label{p3.2}
Let $v$ be a place where both $\eta_i$ are unramified, and
$\mu:F_v^* \to \C^*$ a continuous character with conductor
$\mathfrak{P}^r_v$, $r>0$. If $\Psi^0 _{\eta_v \mu}$ is the
spherical vector and  $\Psi^{\rm new} _{\eta_v \mu}$ the newvector
in $V_{\eta_v \mu}$, then we have

\begin{eqnarray*} \sum_{x \in (\mathcal{O}_v/\mathfrak{P}^r_v)^*}
\mu^{-1}(x) \Psi^0 _{\eta_v \mu}(g
\begin{pmatrix}1&\frac{x}{\pi^r_v}
\\0 & 1\end{pmatrix})= \mu^{-1}(-1)
\frac{\eta_2}{\eta_1}(\mathfrak{P}_v^{r}) \cdot
L_v^{-1}(\frac{\eta_1}{\eta_2},0) \cdot \Psi^{\rm new} _{\eta_v
\mu}(g).\end{eqnarray*}
\end{lem}
\begin{proof}
Put $\Psi''(g):=\sum_{x \in (\mathcal{O}_v/\mathfrak{P}^r_v)^*}
\mu^{-1}(x) \Psi^0 _{\eta_v \mu}(g
\begin{pmatrix}1&\frac{x}{\pi^r_v} \\0 & 1\end{pmatrix})$. To
simplify notation we will write $q$ for $\pi_v^r$ . It is easy to
check that
$$\Psi'' ( \begin{pmatrix}a&b\\0&d
\end{pmatrix} g)=\eta_1(a) \mu(a) \eta_2(d) \mu(d) \Psi''(g)$$ and
we refer to \cite{Be} Lemma 4.4 for the proof of right invariance
under $K^1((q^2))$. This implies that $\Psi'' \in V_{\eta_v
\mu}^{K^1((q^2))}$ is a multiple of the newvector, which is
nonzero only on $B_0(F_v) \begin{pmatrix}
1&0\\\pi_v^r&1\end{pmatrix} K^1((q^2))$. We now give the
calculation of this multiple.

The non-trivial character $\mu$ on $\mathcal{O}_v^*$ descends to a
non-trivial character on $(\mathcal{O}_v/\mathfrak{P}_v^r)^*$,
which implies that $\sum_{x \in
(\mathcal{O}_v/\mathfrak{P}_v^r)^*} \mu^{-1}(x)=0.$ In fact,
$(\mathcal{O}_v/\mathfrak{P}_v^r)^*$ can be replaced by the
subgroups $(1+\mathfrak{P}_v^n)/(1+\mathfrak{P}_v^r)$ for $n=1,
\ldots r-1$ if $r>1$.

We have the Iwasawa decomposition
$$\begin{pmatrix}1&0\\q&1\end{pmatrix}
\begin{pmatrix}1&\frac{x}{q} \\0 & 1\end{pmatrix}=
\begin{pmatrix}1&\frac{x}{q} \\0 & 1\end{pmatrix}
\begin{pmatrix}\frac{1}{1+x} & -\frac{x^2}{q}\\ 0&1+x \end{pmatrix}
\begin{pmatrix}1 & 0\\ \frac{q}{1+x}&1 \end{pmatrix}.
$$
(This works for all $x$ in our sum if we avoid $x =-1$ in our
choice of representatives for $x \in
(\mathcal{O}_v/\mathfrak{P}^r_v)^*$). We obtain
$$\Psi''(\begin{pmatrix} 1&0\\\pi_v^r&1\end{pmatrix})=\sum_{x \in
(\mathcal{O}_v/\mathfrak{P}^r_v)^*}  \mu^{-1}(x)
(\eta_2/\eta_1)(x+1).$$

For $r=1$ we have
$$\Psi''(\begin{pmatrix}
1&0\\\pi_v^r&1\end{pmatrix})=\sum_{x \in
(\mathcal{O}_v/\mathfrak{P}_v)^*, (x+1) \neq \mathfrak{P}_v}
\mu^{-1}(x) + \mu^{-1}(\pi_v-1) (\eta_2/\eta_1)(\pi_v).$$ Using
that $\sum_{x \in (\mathcal{O}_v/\mathfrak{P}_v)^*} \mu^{-1}(x)=0$
this equals
$$ \mu^{-1}(\pi_v-1) \left((\eta_2/\eta_1)(\pi_v)-1 \right)=  \mu^{-1}(-1)
(\eta_2/\eta_1)(\mathfrak{P}_v) \cdot L_v^{-1}(\eta_1/\eta_2,0).$$

For $r>1$ we have
$$\Psi^0 _{\eta_v \mu}(\begin{pmatrix}
1&0\\\pi_v^r&1\end{pmatrix})= \sum_{x \in
(\mathcal{O}_v/\mathfrak{P}^r_v)^*, (x+1) \notin \mathfrak{P}_v}
\mu^{-1}(x) + \sum_{x \in \mathcal{O}_v/\mathfrak{P}^{r-1}_v}
 \mu^{-1}(x \pi_v -1 ) (\eta_2/\eta_1) (x \pi_v).
$$
We rewrite the second sum as
$$ \mu^{-1}(-1)(\eta_2/\eta_1)(\pi_v^r)+\mu^{-1}(-1)
\sum_{n=1}^{r-1} (\eta_2/\eta_1) (\pi_v^{n})
 \sum_{u \in (\mathcal{O}_v/\mathfrak{P}^{r-n}_v)^*}
 \mu^{-1}(1+\pi_v^{n}u).$$
Let $S_n:= \sum_{u \in (\mathcal{O}_v/\mathfrak{P}^{r-n}_v)^*}
 \mu^{-1}(1+\pi_v^{n}u)$. Now we make the following observation: Since $$\sum_{y \in 1+
\mathfrak{P}_v^{m}} \mu^{-1}(y)=0$$ for $m=1, \ldots, r-1$ by our
initial remark we have
$$S_n=
\sum_{w \in \mathcal{O}_v/\mathfrak{P}_v^{r-n}} \mu^{-1} (1+
\pi_v^n w) -
 \sum_{w \in \mathcal{O}_v/\mathfrak{P}_v^{r-n-1}}
 \mu^{-1}(1+\pi_v^{n+1}w)=
  \begin{cases}
    0 & \text{if } n \leq r-2, \\
    -1 & \text{if } n=r-1.
  \end{cases}
 $$
We are left to evaluate
\begin{eqnarray*}\Psi''(\begin{pmatrix}
1&0\\\pi_v^r&1\end{pmatrix})&=& \left( \sum_{x \in
(\mathcal{O}_v/\mathfrak{P}^r_v)^*, (x+1) \notin \mathfrak{P}_v}
\mu^{-1}(x) \right)\\ &+& \mu^{-1}(-1)(\eta_2/\eta_1)(\pi_v^r) -
\mu^{-1}(-1) (\eta_2/\eta_1)(\pi_v^{r-1}).
\end{eqnarray*}
The sum in brackets turns out to be zero as well, since
$$0=\sum_{x \not
\equiv -1 \mod{\mathfrak{P}_v}} \mu^{-1}(x) + \mu^{-1}(-1)
\sum_{n=1}^{r-1} S_n + \mu^{-1}(-1)= \sum_{x \not \equiv -1
\mod{\mathfrak{P}_v}} \mu^{-1}(x).$$ We conclude that
$$\Psi''(\begin{pmatrix}
1&0\\\pi_v^r&1\end{pmatrix})= \mu^{-1}(-1)(\eta_2/\eta_1)(\pi_v^r)
\cdot (1- (\eta_1/\eta_2)(\pi_v)),$$ as desired.
\end{proof}

To translate from the spherical vector to the newvector we
therefore also define the finite twisted sum \begin{equation}
\Psi^{\rm twist}_{\phi_f |\alpha|_f^{z/2}}:= \sum_{v \in S}
\sum_{x \in (\mathcal{O}_v/\mathfrak{P}_v^{r_v})^*} \phi_1^{-1}(x)
\Psi_{\phi_f |\alpha|_f^{z/2}}^0(g
\begin{pmatrix}1&\frac{x}{\pi^{r_v}_v} \\0 & 1\end{pmatrix}_v) \in V_{\phi_f |\alpha|_f^{z/2}}^{K_f^{\rm
new}},\end{equation} where $\mathfrak{P}_v^{r_v} \parallel
\mathfrak{M}_1$. Lemma \ref{p3.2} shows that $\Psi^{\rm
twist}_{\phi_f |\alpha|_f^{z/2}}$ equals $\Psi_{\phi_f
|\alpha|_f^{z/2}}^{\rm new}$ up to essentially the $L$-factors for
$v\in S$.

\subsection{Hecke algebra action} The definition of our special vectors was chosen such
that the Eisenstein cohomology class is a Hecke eigenvector for
almost all $T_{\pi_v}=[K_f^S \begin{pmatrix} \pi_v &0\\0&1
\end{pmatrix} K_f^S]$ (see Section \ref{Scohom} for the definition
of the Hecke operators on relative Lie algebra cocycles). By the
definition of the Eisenstein cohomology class $[{\rm Eis} (\phi,
\Psi^0_{\phi})]$, it suffices to check the effect of the Hecke
operator $T_{\pi_v}$ on $\Psi^0_{\phi_f} \in V_{\phi_f}^{K_f^S}$.
We note that for any $\Psi= \prod_w \in V_{\phi_f}^{K_f}$ with
$K_f=\prod_w K_w$ the action of $T_{\pi_v}$ is described by
$(T_{\pi_v}.\Psi)(g)= \prod_{w \neq v} \Psi_w(g_w) \cdot \left(
\sum_{i } \Psi_v(g_v \gamma_i) \right),$ where $K_v
\begin{pmatrix} \pi_v &0\\0&1 \end{pmatrix} K_v= \coprod_i \gamma_i
K_v$. By \cite{Cas} we know that
$V_{\phi_v}^{K^1(\mathfrak{P}_v^s)}$ for $\mathfrak{P}_v^s \|
\mathfrak{M}_1 \mathfrak{M}_2$ is 1-dimensional so we get for $v
\notin S$ that
$$T_{\pi_v} (\Psi^0_{\phi_f})=a_v(\phi) \Psi^0_{\phi_f} \text{ for some }
a_v(\phi) \in \C.$$
\begin{lem} \label{Hecke}
For $v \notin S$ the class $[{\rm Eis} (\phi, \Psi^0_{\phi_f})]$
is an eigenvector for $T_{\pi_v}$.
\end{lem}
 If in addition $v \nmid \mathfrak{M}$, for example,  the eigenvalue is
\begin{eqnarray*} T_{\pi_v}(\Psi^0_{\phi_f})(1)&=& \Psi^{\rm
new}_v(\begin{pmatrix}1&0\\0&\pi_v \end{pmatrix}) \ +\sum_{a\in
\mathcal{O}_v/\mathfrak{P}_v} \Psi^{\rm
new}_v(\begin{pmatrix}\pi_v & a\\0&1\end{pmatrix})
\\&=&\phi_{2,v}(\mathfrak{P}_v)+ {\rm Nm}(\mathfrak{P}_v)
\phi_{1,v}(\mathfrak{P}_v).\end{eqnarray*} For the calculation of
the eigenvalues in other cases see \cite{Be}  Lemma 3.11.

\subsection{Constant terms}
We will be interested in the image of the Eisenstein cohomology
classes under the restriction map
$${\rm res}: H^1(S_{K_f},
\widetilde{M_{\C}})\cong H^1(\overline S_{K_f},
\widetilde{M_{\C}})\to H^1(\partial \overline S_{K_f},
\widetilde{M_{\C}}) \cong H^1(\partial \tilde S_{K_f},
\widetilde{M_{\C}}).$$ We are in the case that Harder calls
``balanced" where the restriction map maps diagonally into the
cohomology of the boundary, in the sense that on the level of
functions
$${\rm res} \circ {\rm Eis}: V_{\phi}^{K_f} \to V_{\phi}^{K_f}
\oplus V_{w_0.\phi}^{K_f}$$
$$\Psi \mapsto \Psi + \star \frac{L(-1,\phi_1/
\phi_2)}{L(0, \phi_1/\phi_2)}T_{\phi} \Psi \in V_{\phi}^{K_f}
\oplus V_{w_0.\phi}^{K_f},$$ for an intertwining operator
$T_{\phi} :V_{\phi}^{K_f} \to V_{w_0.\phi}^{K_f}$ and some
non-zero factor $\star$.

By \cite{Ha79} Proposition 1.6.1, or \cite{Z} Proposition 2.2.3
(see also \cite{Schw} Satz 1.10 in the case of automorphic forms)
the restriction of a cohomology class in $H^1(S_{K_f},
\widetilde{M_{\C}})$, represented by a relative Lie algebra
1-cocycle $$\omega \in {\rm
Hom}_{K_{\infty}}(\mathfrak{g}_{\infty}/\mathfrak{k}_{\infty},
 C^{\infty} (G(\Q) \backslash G(\A)/K_f) \otimes M_{\C}),$$
is given by the class of the constant term $\omega_B$, where the
\textit{constant term with respect to a parabolic} $P$ is defined
as
$$\omega_P(g)=\int_{U_P(\Q) \backslash U_P(\A)} \omega(ug) \, du$$
for an appropriate Haar measure $du$. Recall the decomposition of
$$\partial \tilde S_{K_f}=B(\Q) \backslash G(\A)/K_f K_{\infty}$$
into its connected components given in (\ref{BSbdry}). For the
parabolic $B^{\eta}$ with $\eta \in G(\Q)$ and $\gamma \in
G(\A_f)$ let
$${\rm res}_{B^{\eta}}^{\gamma}: H^1(S_{K_f},
\widetilde{M_{\C}}) \to  H^1(\Gamma_{\gamma, B^{\eta}} \backslash
\h, \widetilde{M_{\C}})$$ be the restriction map to the boundary
component $\Gamma_{\gamma, B^{\eta}} \backslash \h
\overset{j_{\eta, \gamma}}{\hookrightarrow}
\partial \tilde S_{K_f}$. It is easy to check that $${\rm res}_{B^{\eta}}^{\gamma}
[\omega]=j^*_{\gamma} [\omega_{B^{\eta}}]=j^*_{\eta, \gamma}
[\omega_B].$$  It suffices therefore to calculate the constant
term $\omega_B$.
\begin{prop} \label{const}

  \emph{(a)} If $\Psi=\Psi^0_{\phi_f |\alpha|_f^{z/2}}$ then
$${\rm Eis}(\Psi )_B= \omega_z(\Psi)+ c(\phi, z) \omega_{-z}(\Psi^0_{w_0.(\phi_f
|\alpha|_f^{z/2})}),$$ where $$c(\phi,z)=(d_F)^{-1/2} \frac{2
\pi}{z+m+1} (-1)^{n+1} \cdot \frac{L(z-1, \chi)}{L(z, \chi)} \cdot
 \prod_{v | \mathfrak{M}} c_v(\phi,z)$$ with $c_v(\phi,z):=
 \int_{U_0(F_v)} \Psi^{\rm new}_v(w_0 u_v\begin{pmatrix} 1 & 0 \\ \pi_v^{t_v} & 1
\end{pmatrix}) du_v$, where $\mathfrak{P}_v^{t_v} \parallel \mathfrak{M}_2$.
If only one of $\{\phi_{1,v},\phi_{2,v}\}$ is
ramified then $c_v(\phi,z)= \frac{\phi_{2,v}(-1)}{{\rm
Nm}(\mathfrak{M}_{1,v})}$.

  \emph{(b)} If $\Psi=\Psi^{\rm new}_{\phi_f |\alpha|_f^{z/2}}$ then
$${\rm Eis}(\Psi )_B= \omega_z(\Psi)+ c'(\phi, z) \omega_{-z}(\Psi^{\rm new}_{w_0.(\phi_f
|\alpha|_f^{z/2})}),$$ where
$$c'(\phi,z)=c(\phi,z) \cdot \prod_{v \in S} [ (1-{\rm
Nm}(\mathfrak{P}_v)) \cdot (\chi/{\rm Nm})(\mathfrak{P}_v^{r_v})
\cdot L_v(z, \chi) ],$$ with $\mathfrak{P}_v^{r_v} \parallel
\mathfrak{M}_1$.
\end{prop}
\begin{proof}
See \cite{HaGL2} Theorem 2(3) and \cite{Be} Proposition 3.5.
\end{proof}

\subsection{Translation to group cohomology}
The decomposition (\ref{BSbdry}) and Proposition \ref{dRgroup}
imply that we have an isomorphism
\begin{equation}\label{bdry2gpcoh} H^1(\partial \tilde S_{K_f}, \widetilde M_{\C}) \cong
\bigoplus_{[{\rm det}(\gamma)] \in \pi_0(K_f)} \bigoplus_{[\eta]
\in \mathbf{P}^1(F)/\Gamma_{\gamma}} H^1(\Gamma_{\gamma,
B^{\eta}}, M_{\gamma} \otimes \C ).\end{equation} The following
Lemma calculates the image of the boundary cohomology class we
defined in Section \ref{s3.1} under this isomorphism:
\begin{lem} \label{omega0_group}
 For $\phi: T(\Q)
\backslash T(\A) \to \C^*$, $\Psi \in V_{\phi_f}^{K_f}$, and
$\omega=\omega_0(\phi,\Psi)$ as defined in (\ref{omega0}) the
image of $[\omega]$ in $H^1(\Gamma_{\gamma, B^{\eta}}, M_{\C})$ is
represented  by the 1-cocycle
 $$\eta^{-1}_{\infty}
\begin{pmatrix}a&x\\0&d
\end{pmatrix} \eta_{\infty} \mapsto  \Psi(\eta_f \gamma)
  \begin{cases}
    x \int_0^1 \begin{pmatrix} 1 & tx\\
0&1\end{pmatrix}.Y^m \overline X^n dt& \text{ if } \phi \text{ of infinity type (A)}, \\
    \overline x \int_0^1 \begin{pmatrix} 1 & tx\\
0&1\end{pmatrix}.X^m \overline Y^n dt & \text{ if } \phi \text{ of
infinity type (B)}.
  \end{cases}$$
(Here we denote by $\eta_f$ and $\eta_{\infty}$ the images of
$\eta \in G(\Q)$ in $G(\A_f)$ and $G_{\infty}$, respectively.)
\end{lem}

\begin{proof}
Put $P=B^{\eta}$. Recall from Proposition \ref{dRgroup} that for
the basepoint $x_0=\eta_{\infty}^{-1} K_{\infty}$ the image of
$[\omega]$ is represented by the cocycle given on $U^{\eta}$ by
$$\mathcal{G}_{x_0}
( j_{\eta,\gamma}^*(\tilde \omega)) (\eta^{-1}_{\infty}
\begin{pmatrix}1&x\\0&1
\end{pmatrix} \eta_{\infty} )= \int_{\eta^{-1}_{\infty} K_{\infty}}^{\eta^{-1}_{\infty}
\left(\begin{smallmatrix}1&x\\0&1 \end{smallmatrix} \right
)K_{\infty}}
j_{\eta,\gamma}^*(\tilde\omega)=\int_{(K_{\infty},\eta_f
\gamma)}^{ (\left(\begin{smallmatrix}1&x\\0&1
\end{smallmatrix} \right )K_{\infty},\eta_f \gamma)}
\tilde\omega.$$ Here $\tilde \omega \in \Omega^1(G(\A)/K_f
K_{\infty}\otimes M_{\C})^{B(\Q)}$ is the closed 1-form associated
to $\omega$, given by
$$\tilde
\omega(x)(T):=g_{\infty}.\omega(g_{\infty}, g_f)({\rm D}_{{\rm
L}_{g_{\infty}}^{-1}} T)$$ for $x=(x_{\infty}, g_f)=(g_{\infty}
K_{\infty}, g_f) \in \h \times G(\A_f)/K_f$ and $T \in
T_{x_{\infty}} \h$, independent of the choice of $g_{\infty}$. To
calculate the path integral we apply the following lemma, adapted
from \cite{Wes}:

\begin{lem}[(\cite{Wes} Lemma 5.1)]
Given $h:\Real \to G_{\infty}$ a differentiable homomorphism and
$g_{\infty}\in G_{\infty}$, $g_f\in G(\A_f)/K_f$, define $c: \Real
\to \h$ by $c(t):=h(t) \cdot g_{\infty} \cdot K_{\infty}$. For
$a_0, a_1 \in \Real$ let $y_i:=c(a_i)$ and denote $\dot{h}:=(Dh)_0
T_0 \in \mathfrak{g}_{\infty}$. Then one has the following
equality:
$$\int_{y_0}^{y_1} \tilde \omega = \int_{a_0}^{a_1}
(h(t)g_{\infty}).\omega(h(t)g_{\infty},
g_f)(g^{-1}_{\infty}\dot{h} g_{\infty}) dt. $$
\end{lem}
We take $y_0=K_{\infty}$, $g_{\infty}=1$, $g_f=\eta_f\gamma$,
$h(t)=
\begin{pmatrix}1&xt \\ 0 &1
\end{pmatrix} \in G_{\infty}$, $a_0=0$, $a_1=1$, and
obtain:
$$\mathcal{G}_{x_0}
(j_{\eta,\gamma}^*(\tilde \omega)) (\eta^{-1}_{\infty}
\begin{pmatrix}1&x\\0&1
\end{pmatrix} \eta_{\infty} )=
\int_0^1 (\begin{pmatrix} 1 & tx\\
0&1\end{pmatrix}).\omega(\begin{pmatrix} 1 & tx\\
0&1\end{pmatrix}, \eta_f \gamma) ( \begin{pmatrix}
0&x\\0&0\end{pmatrix}) dt.$$ One calculates that for $x \in \C$
$$x S_+ -\overline x S_{-}= \begin{pmatrix} 0&x/2\\\overline x/2&
0
\end{pmatrix} \equiv \begin{pmatrix} 0&x\\0&0\end{pmatrix}
\mod{\mathfrak{k}_{\infty}}.$$ To complete the proof one checks
that $\mathcal{G}_{x_0}(j_{\eta,\gamma}^*(\tilde \omega))$ is
always zero on $\eta_{\infty}^{-1}
\begin{pmatrix} a & 0 \\ 0 & d\end{pmatrix} \eta_{\infty}$ since $\omega$
vanishes along $H \in \mathfrak{p}_{\C}$.
\end{proof}

\subsection{Rationality of Eisenstein cohomology class and integrality of constant term}
Harder proves that the transcendentally defined Eisenstein
operator is, in fact, rational, i.e., that ${\rm Eis}$ is defined
over $\overline \Q$ so that we have
$$V_{\phi_f}^{K_f} \to H^1(S_{K_f}, \widetilde M_{\overline \Q})$$
$$ \Psi \mapsto [{\rm Eis}(\phi, \Psi)]$$ and for $\sigma \in {\rm
Gal}(\overline \Q/\Q)$ and $\Psi \in V_{\phi_f}^{K_f}$ the
equation $[{\rm Eis}(\phi, \Psi)]^{\sigma}=[{\rm
Eis}(\phi^{\sigma}, \Psi^{\sigma})]$ holds (see \cite{HaGL2}
Corollary 4.2.1(a); the proof uses the ``Multiplicity one Theorem"
for automorphic forms and the vanishing of residual interior
cohomology in our case).

We are interested in the $p$-adic properties of the Eisenstein
cohomology class. From now on, let $\phi=(\phi_1,
\phi_2):T(\Q)\backslash T(\A) \to C^*$ denote a continuous
character of infinity type (A) for which the conductors of both
$\phi_i$ are coprime to $(p)$, and let $m \geq n$ (recall that
$\chi:=\phi_1/\phi_2$). Let $\mathcal{O}_{\phi}$ denote the ring
of integers in the finite extension $F_{\phi}$ of
$F_{\mathfrak{p}}$ obtained by adjoining the values of the finite
part of both $\phi_i$ and $L^{\mathrm{alg}}(0,\chi)$. For the
definition of the latter and its $p$-adic properties see Theorem
\ref{Linteg}. Then the above discussion shows:
\begin{prop} \label{ration} For $\Psi=\Psi^0_{\phi_f}$ or $\Psi^{\rm
new}_{\phi_f}$ we have
$$[{\rm Eis}(\phi, \Psi)] \in
H^1(S_{K_f},\widetilde M_{F_{\phi}})$$ for $K_f=K_f^S$ or
$K_f=K_f^{\rm new}$, respectively.
\end{prop}

\begin{defn} \label{denom}
If $\mathcal{O}_L$ is the ring of integers in a local field $L$
over $F_{\mathfrak{p}}$, define for any $c \in H^1(
S_{K_f},\widetilde M_L)$ the \textit{denominator (ideal)}
$$\delta(c)=\{a\ \in \mathcal{O}_L: ac \in H^1(
S_{K_f},\widetilde M_{\mathcal{O}_L})_{\rm free}  \}.$$ Here we
identify $H^1( S_{K_f}, \widetilde M_{\mathcal{O}_L})_{\rm free}$
with $\mathrm{im} (H^1( S_{K_f},\widetilde M_{\mathcal{O}_L}) \to
H^1( S_{K_f},\widetilde M_L))$.
\end{defn}

In this paper we prove (under certain conditions on the conductors
of the characters $\phi_i$) that $L^{\rm alg}(0,\chi) \cdot
\mathcal{O}_{\phi}$ is a lower bound on the denominator of the
Eisenstein cohomology class.

For the arithmetic of the Eisenstein cohomology class its constant
term is, of course, of great importance. Put
$N:=M(m,n,-m-k,-n-\ell)$. Note that since $N^{\vee} \cong M=M(m,n,
k, \ell)$ as $G(\Q)$-modules, $[{\rm Eis}(\phi, \Psi^0_{\phi_f})]$
can also be interpreted as a cohomology class for the local system
$\widetilde N^{\vee}_{\C}$. The following result shows that under
certain conditions its constant term is integral already with
respect to the lattice $(N_{\mathcal{O}})^{\vee} \subset
M_{\mathcal{O}} \subset M$. We first prove the following Lemma:
\begin{lem}
Assume in addition that $p>m+1$. We have
$$[\omega_0(\phi, \Psi_{\phi}^0)], [\omega_0(w_0.\phi,
\Psi_{w_0.\phi}^0)]\in H^1(\partial \tilde S_{K_f^S},
\widetilde{(N_{\mathcal{O}_\phi})^{\vee}})_{\rm free}.$$
\end{lem}
\begin{proof}
We recall from Section \ref{Scohom} that we have an $R$-functorial
isomorphism
$$H^1(\partial \tilde S_{K_f^S}, \widetilde N^{\vee}_{R}) \cong
\bigoplus_{[\gamma] \in \pi_0(K_f^S)}\bigoplus_{[\eta] \in
\mathbf{P}^1(F)/\Gamma_{\gamma}} H^1(\Gamma_{\gamma, B^{\eta}},
(N_{\gamma})^{\vee}\otimes  R).$$ We will show that for
$[\omega_0(\phi, \Psi_{\phi}^0)] \in H^1(\partial \tilde
S_{K_f^S}, \widetilde N^{\vee}_{\C})$ the restrictions to each of
the summands on the right hand side lie in the image of
$H^1(\Gamma_{\gamma, B^{\eta}}, (N_{\gamma})^{\vee} \otimes
\mathcal{O}_{\phi})$ inside $H^1(\Gamma_{\gamma, B^{\eta}} ,
(N_{\gamma})^{\vee} \otimes {\C})$. We showed in Lemma
\ref{omega0_group} that for each $\gamma$ and $\eta$ this
restriction is given by
$$\eta^{-1} \begin{pmatrix} a&x\\0 &b\end{pmatrix} \eta \mapsto x
\Psi^0_{\phi_f}(\eta_f \gamma) \int_0^1
\begin{pmatrix} 1&tx\\0 &1\end{pmatrix}.Y^m \overline X^n dt.$$ To
check the integrality we choose representatives $\gamma$ and
$\eta$ whose $\mathfrak{p}$- and
$\overline{\mathfrak{p}}$-components are units (i.e., they are
elements of ${\rm GL}_2(\mathcal{O}_p):=\prod_{v \mid p} {\rm
GL}_2(\mathcal{O}_v)$). This is possible for $\gamma$ by the
Chebotarev density theorem. For $\eta$ this follows from
$${\rm GL}_2(F_p):=\prod_{v \mid p} {\rm GL}_2(F_{v}) =B_0(F) {\rm
GL}_2(\mathcal{O}_p).$$

For such $\gamma$ we get $N_{\gamma}^{\vee} \otimes
\mathcal{O}_{\phi}=(N_{\mathcal{O}_{\phi}})^{\vee}$ and we need to
show that the values of the group cocycle satisfy that the
coefficient of $X^iY^{m-i}\overline X^j \overline Y^{n-j}$ lies in
$\begin{pmatrix}
    m \\
    i
  \end{pmatrix}
  \begin{pmatrix}
    n \\
    j
  \end{pmatrix} \mathcal{O}_{\phi}$. Firstly,
  $\Psi^0_{\phi_f}(\eta_f \gamma) \in \mathcal{O}_{\phi}^*$ by the
  definition of the spherical vector at
places away from the conductors of the $\phi_i$ together with
Lemma \ref{charinteg}. We also note that $\Gamma_{\gamma} \cap
G(\Q_p) \subset {\rm GL}_2(\mathcal{O}_p)$, so $x$ lies in
$\mathcal{O}_{\phi}$. Lastly, the integral provides us with the
correct coefficients for the monomials up to $p$-adic units if we
assume $p>m+1$.

A similar argument for $[\omega_0(w_0.\phi, \Psi_{w_0.\phi}^0)]$
proves integrality if $p>n+1$ (and $m\geq n$ by assumption).
\end{proof}

We proved in Proposition \ref{const} that ${\rm
Eis}(\Psi^0_{\phi_f})_B=\omega_0(\Psi^0_{\phi_f}) + c(\phi, 0)
\omega_0(\Psi^0_{w_0.\phi_f}).$ By \cite{HaGL2} Corollary 4.2.2 we
know that $c(\phi,0) \in F_{\phi}$. The following proposition
analyzes conditions when $c(\phi,0)$ (and so by the preceding
Lemma the constant term of the Eisenstein class) is integral.
\begin{prop} \label{resinteg}
Assume in addition that $p>m+1$ and that the conductors of
$\phi_1$ and $\chi=\phi_1/\phi_2$ are coprime.

  \emph{(a)}
The constant term $[{\rm Eis}(\phi, \Psi^0_{\phi_f})_B]$ is
integral if and only if $\frac{L^{\rm alg}(-1, \chi)}{L^{\rm
alg}(0, \chi)} \in \mathcal{O}_{\phi}$.

    \emph{(b)}
For $m=n$, we have $[{\rm Eis}(\phi, \Psi^0_{\phi_f})_B] \in
H^1(\partial \overline S_{K_f^S},
\widetilde{(N_{\mathcal{O}_\phi})^{\vee}})_{\rm free}$ if and only
if $\frac{L(0, \overline \chi)}{L(0,\chi)}\in \mathcal{O}_{\phi}$.

    \emph{(c)}
If $m=n$ and $\chi^c(x):=\chi(\overline x)=\overline \chi(x)$ for
all $x \in \A_F^*$ then
$$[{\rm Eis}(\phi, \Psi^0_{\phi_f})_B] \in H^1(\partial \overline
S_{K_f^S}, \widetilde{(N_{\mathcal{O}_\phi})^{\vee}})_{\rm
free}.$$
\end{prop}

\begin{rem}
\begin{enumerate}
  \item In \cite{Be} we called characters $\chi$ satisfying
$\chi^c=\overline \chi$ \textit{anticyclotomic}. These include, in
particular, the everywhere unramified characters.

  \item  As explained in the introduction, the integrality of the constant term of the Eisenstein cohomology
class (together with the bound on the denominator) is interesting
for finding congruences between cuspidal and Eisenstein cohomology
classes controlled by the $L$-value $L^{\rm alg}(0,\chi)$. Since
cuspidal cohomology classes only exist for $m=n$ by Wigner's Lemma
(see, e.g. \cite{Ha06} \S 3) these are the coefficient systems we
are most interested in.
\item By considering $c(\phi,0) \overline{c(\phi,0)}$ and using $\overline{c(\phi,0)}=c(\overline
\phi,0)$ (cf. \cite{HaGL2} Corollary 4.2.2) one can show that the
quotient of $L$-values in (b) is always integral for either
$\mathfrak{p}$ or $\overline{\mathfrak{p}}$.
\end{enumerate}
\end{rem}

\begin{proof}
One easily checks that
$$c(\phi,0)=\frac{L^{\rm alg}(0, \chi | \cdot |^{-1})}{L^{\rm alg}(0,
\chi)} (-1)^{n+1} \prod_{v \mid \mathfrak{M}} c_v(\phi,0).$$ Since
$\mathfrak{M}$ is coprime to $p$ and the conductor of $\phi_1$ we
have $c_v(\phi,0) =\frac{\phi_{2,v}(-1)}{{\rm
Nm}(\mathfrak{M}_{1,v})} \in \mathcal{O}_{\phi}^*$.

For (b) we can apply the functional equation (cf. \cite{dS} p. 37)
and use that $m=n$, and therefore $\chi \overline \chi= | \cdot
|^2$, to obtain $$c(\phi,0)= \frac{L(0, \overline
\chi)}{L(0,\chi)} (-1)^{n+1} W(\chi) \sqrt{{\rm Nm}(\mathfrak{M})}
\prod_{v \mid \mathfrak{M}} c_v(\phi,0),$$ where $W(\chi)$ is the
Artin root number for $\chi$ (see Section \ref{Hcharacters}). By
the assumption on the conductor of $\chi$ both $\sqrt{{\rm Nm}(
\mathfrak{M})} $ and the root number $W(\chi)$ lie in
$\mathcal{O}_{\phi}^*$.

Lastly, if $\chi^c=\overline \chi$ then we have $L(0,\overline
\chi)=L(0, \chi^c)=L(0,\chi)$, so $\frac{L(0, \overline
\chi)}{L(0,\chi)}=1$.

\end{proof}

\section{Toroidal integral}
In this section we calculate the integral of twists of the
Eisenstein cocycle defined in Section 3 against certain relative
cycles.
\subsection{Definition of relative cycles} \label{3.2.1}
Recall that strong approximation implies that $S_{K_f}$ is the
finite disjoint union of connected components indexed by a set of
representative $\{[\xi]\}$ for $\pi_0(K_f)$ with $\xi\in
\A_{F,f}^*$. The connected component lying above $\xi$ is given by
$\Gamma_{\xi} \backslash \h$ with $\Gamma_{\xi}:= G(\Q) \cap
\begin{pmatrix} 1&0\\0&\xi \end{pmatrix} K_f \begin{pmatrix}
1&0\\0&\xi \end{pmatrix}^{-1}$ and we defined in Section
\ref{s2.2} an embedding $j_{\xi}:=j_{\begin{pmatrix} 1&0\\0&\xi
\end{pmatrix}}: \Gamma_{\xi} \backslash
\h \hookrightarrow S_{K_f}$.  For each $\xi \in \A_{F,f}^*$ let
$$\sigma_{\xi}=j_{\xi} \circ \tau : \C^* \to S_{K_f},$$ where
$\tau: \C^* \to G_{\infty}: z \mapsto
\begin{pmatrix}1&0 \\0 & z \end{pmatrix}$. We will use $\sigma_{\xi}$ to denote both
the map to $G(\A)$ and the induced map to $S_{K_f}$. We consider
the path in $S_{K_f}$ given by $\sigma_{\xi}|_{\Real_{>0}^*}$,
which is the restriction of a path (also denoted by
$\sigma_{\xi}$) in $\overline S_{K_f}$: for the component
$\Gamma_{\xi} \backslash \oh$ that path is $\sigma_{\xi}: [0,
\infty] \to \Gamma_{\xi} \backslash \oh \subset \overline S_{K_f}:
t\mapsto j_{\xi}((t,0)).$
 The paths $\sigma_{\xi}$ are not
1-cycles in $\overline S_{K_f}$. They are, however, relative
cycles in $C_1( \Gamma_{\xi} \backslash \oh ,\partial (
\Gamma_{\xi} \backslash \oh), \Z)$ (cf. \cite{Ko} \S 5.2). Since
the endpoints lie in the $\infty$- and $0$-cusps (use
$\begin{pmatrix} 1&0\\0&\frac{1}{s}
\end{pmatrix}K_{\infty}=w_0 \begin{pmatrix} 1&0\\0&s
\end{pmatrix} K_{\infty}$) they are, in fact, relative cycles for $$H_1(\Gamma_{\xi}
\backslash \oh, \Gamma_{\xi, B} \backslash e(B) \cup \Gamma_{\xi,
B^w} \backslash e(B^w), \Z).$$

Put $$Q_{m',n'}:=X^{m-m'} Y^{m'} \overline X^{n-n'} \overline
Y^{n'} \in N_{\mathcal{O}} \text{ for } 0 \leq m' \leq m, \, 0
\leq n' \leq n.$$ Let \begin{equation} \label{theta}
\theta_{m',n'}:F^* \backslash \A_F^* \to \C^* \end{equation} be a
Hecke character with $\theta_{m',n',\infty}(z)=z^{m-m'+k}
\overline z^{n-n'+\ell}$ and conductor $\mathfrak{N}$ coprime to
$(p d_F)$ and the conductors of $\phi_1$ and $\phi_2$ such that
$\#(\mathcal{O}/\mathfrak{N})^*$ is also coprime to $p$. Denote by
$T$ the set of places where $\theta_{m',n'}$ is ramified.

For any $\xi \in \A_{F,f}^*$ consider now the chain
$$\left( \sigma_{\xi} \otimes Q_{m',n'} \right) \cdot \theta_{m',n'}
(\xi) \in C_1(\Gamma_{\xi} \backslash \h, N_{\xi}\otimes
\mathcal{O}_{\theta}),$$ with $N_{\xi}:=N_{\begin{pmatrix}
1&0\\0&\xi
\end{pmatrix}}\cong j_{\xi}^* \widetilde N_{\mathcal{O}}$ and $\mathcal{O}_{\theta}$
the ring of integers in the finite extension $F_{\theta}$ of
$F_{\phi}$ obtained by adjoining the values of $\theta_{m',n',f}$.
Here we use Lemma \ref{charinteg} to check the integrality of the
chain. Since chains for $S_{K_f}$ are defined as
$G(\Q)$-coinvariants of chains in $G(\A)/K_f K_{\infty}$ the sum
$$\sum_{[\xi] \in \pi_0(K_f)} \left( \sigma_{\xi} \otimes
Q_{m',n'} \right) \cdot \theta_{m',n'} (\xi) \in C_1(S_{K_f},
N_{\mathcal{O}} \otimes \mathcal{O}_{\theta})$$ is independent of
the choice of representatives $\xi$ for the connected components.
As observed above this chain is, in fact, a relative cycle with
endpoints in the $\infty$ and $0$-cusps of all the connected
components $\Gamma_{\xi} \backslash \oh$.

\subsection{Twisted version of Eisenstein cocycle}
We now define the following twisted version of the Eisenstein
cocycle: Let $\eta=(\eta_1, \eta_2): F^*\backslash \A_F^* \to
\C^*$ be a continuous character. For $\Psi \in V_{\eta_f}$ let
$${\rm Eis}^{\theta}(g, \Psi):= \sum_{v \in T} \sum_{x \in
(\mathcal{O}_v/\mathfrak{N}_v)^*} \theta_{m',n'}^{-1}(x) {\rm
Eis}\left(g \begin{pmatrix} 1 & -\frac{x}{\pi_v^{{\rm ord}_v
\mathfrak{N}}} \\ 0 &1
\end{pmatrix}_v, \Psi \right).$$ Note that if $\Psi=\prod_v
\Psi_v$ then the twisting can also be done on the level of the
function, i.e., $${\rm Eis}^{\theta}(g, \Psi)={\rm Eis}(g,
\Psi^{\theta}), $$ where $\Psi^{\theta}(g):=\prod_{v \in T}
\Psi^{\theta}_v(g_v) \prod_{v \notin T} \Psi_v(g_v)$ and
$$\Psi^{\theta}_v(g_v):= \sum_{v \in T} \sum_{x \in
(\mathcal{O}_v/\mathfrak{N}_v)^*} \theta_{m',n'}^{-1}(x)
\Psi_v\left(g \begin{pmatrix} 1 & -\frac{x}{\pi_v^{{\rm ord}_v
\mathfrak{N}}} \\ 0 &1
\end{pmatrix}_v \right).$$
Lemmata \ref{l3.2} and \ref{p3.2} imply:
\begin{lem} \label{thetatwist}
For $v\in T$ we have
$$\Psi^{{\rm new}, \theta}_{\eta, v}(g)=\Psi^{\rm new}_{\eta
\theta, v}(g) \theta_v^{-1}(-{\rm det}(g)) \cdot
(\eta_2/\eta_1)(\pi_v^{{\rm ord}_v \mathfrak{N}}) \cdot
L_v^{-1}(\eta_1/\eta_2,0) .$$
\end{lem} Now consider $\eta=\phi|\alpha|^{z/2}$. We note that if $[{\rm Eis}(g, \Psi)] \in H^1(S_{K_f}, \widetilde
M_{F_{\phi}})$ then
$$[{\rm Eis}^{\theta}(g, \Psi)] \in H^1(S_{(K_f)^{\theta}}, \widetilde M_{F_{\theta}}),$$
where for $K_f=\prod_v K_v$ we define $(K_f)^{\theta}=\prod_v
K_v'$ with $K_v'=K_v$ for $v \notin T$ and $K_v'=K_v \cap
U^1(\mathfrak{N}_v)$ for $v\in T$. We will see in Lemma
\ref{relativecocycle} that this twisting makes ${\rm
Eis}^{\theta}(\Psi)$ a relative cocycle with respect to the $0$-
and $\infty$-cusps of each connected component (in the sense of
Proposition \ref{Borel}).

\subsection{Calculation of the toroidal integral}

We now want to integrate ${\rm Eis}^{\theta}(\Psi)$ for $\Psi \in
V_{\phi_f |\alpha|^{z/2}_f}^{K_f^{\rm new}}$ over the relative
cycle defined in Section \ref{3.2.1}: Put $K_f^{\theta}:=(K_f^{\rm
new})^{\theta}$ and let
\begin{eqnarray*}I(\phi, \theta,\Psi, z)&:=& \sum_{[\xi] \in
\pi_0(K^{\theta}_f)} \theta_{m',n'}(\xi) \int_{\sigma_{\xi}
\otimes Q_{m',n'}} {\rm
Eis}^{\theta}(\Psi)\\
&=&\sum_{[\xi] \in \pi_0(K^{\theta}_f)}\theta_{m',n'}(\xi)
\int_0^{\infty} \left \langle Q_{m',n'}, {\rm Eis}^{\theta}(
\sigma_{\xi}(t),\Psi)(d\sigma_{\xi}(t \frac{\partial}{\partial
t})) \right \rangle \frac{dt}{t}.\end{eqnarray*} We will first
evaluate this ``toroidal integral" for $\Psi=\Psi^{\rm new}:=
\Psi^{\rm new}_{\phi_f |\alpha|_f^{z/2}}$. Rewriting the
Eisenstein cocycle as a relative Lie algebra cocycle we have
\begin{eqnarray*}\lefteqn{I(\phi, \theta, \Psi^{\rm new}, z)=}\\&=\sum_{[\xi] \in
\pi_0(K^{\theta}_f)} \theta_{m',n'}(\xi)  \int_0^{\infty} \left
\langle \begin{pmatrix} 1&0
\\0&t^{-1}\end{pmatrix}. Q_{m',n'}, {\rm Eis}^{\theta}(
\sigma_{\xi}(t),\Psi^{\rm new}
)(d\sigma_{\xi}(\frac{\partial}{\partial t} |_{t=1}))
 \right \rangle \frac{dt}{t}.\end{eqnarray*} Using the $K_{\infty}$-invariance of the
Eisenstein cocycle, the argument on pp.107/8 in \cite{Ko} shows
that this equals
\begin{eqnarray*}\lefteqn{\sum_{[\xi] \in \pi_0(K^{\theta}_f)} \theta_{m',n'}(\xi)
\int_0^{2\pi} \int_0^{\infty}}\\ &\left \langle
\begin{pmatrix} 1&0 \\0&u^{-1}\end{pmatrix}. Q_{m',n'},
 {\rm Eis}^{\theta}(\sigma_{\xi}(u),\Psi^{\rm new} ) ({\rm
Ad}(\sigma_{\xi}(e^{-i \varphi}))
d\sigma_{\xi}(\frac{\partial}{\partial t} |_{t=1})) \right \rangle
\frac{dt}{t} \wedge \frac{d \varphi}{2 \pi}\end{eqnarray*} with
$u=t e^{i\varphi} \in \C^*$. Since
$d\sigma_{\xi}(\frac{\partial}{\partial t} |_{t=1})=\frac{H}{2}$,
this equals
$$\sum_{[\xi]  \in \pi_0(K^{\theta}_f)} \theta_{m',n'}(\xi) \int_{\C^*} \theta_{m',n',\infty}(x_{\infty})
\left \langle Q_{m',n'}, {\rm Eis}^{\theta}( \begin{pmatrix}1&0
\\0 & \xi x_{\infty}\end{pmatrix},\Psi^{\rm new} )(\frac{H}{2})
\right \rangle d^*x_{\infty}$$ with $d^*x_{\infty}:=\frac{i}{4
\pi} \frac{dx_{\infty} \wedge d\overline x_{\infty}}{x_{\infty}
\overline x_{\infty}}$.

Note that ${\rm det}(K_f^{\theta})= \prod_{v\notin T}
\mathcal{O}_v^* \prod_{v\in T} (1 + \mathfrak{N}_v)$. Normalize a
Haar measure $d^*x=d^*x_{\infty} \prod_{v \nmid \infty} d^*x_v$ on
$\A_F^*$ such that for finite places
$\int_{(\mathcal{O}_v/\mathfrak{N}_v)^*} d^*x_v =1$. Using the
right $K_f^{\theta}$-invariance we can then write
$$I(\phi,\theta,\Psi^{\rm new}, z)=\int_{F^* \backslash \A_F^*} \theta_{m',n'} (x) \left \langle Q_{m',n'},
{\rm Eis}^{\theta}(\begin{pmatrix}1&0 \\0 & x \end{pmatrix},
\Psi^{\rm new})(\frac{H}{2}) \right \rangle  d^*x.$$

\begin{prop} \label{p3.4}
For ${\rm Re}(z)\gg 0$ the integral $I(\phi,\theta,\Psi^{\rm new},
z)$ converges and the value is
\begin{eqnarray*} &&\frac{L(\frac{z}{2}, \phi_1 \theta_{m',n'} )
L(\frac{z}{2}, \phi_2^{-1} \theta^{-1}_{m',n'})}{L^{S}(z,
\phi_1/\phi_2)} \cdot
\frac{\Gamma(\frac{z}{2}+m-m'+1)\Gamma(\frac{z}{2}+m'+1)}{\Gamma(z+m+2)}
\cdot \#(\mathcal{O}/\mathfrak{N})^* \cdot
\\ &\cdot& \frac{(-1)^{n-n'+k+\ell}}{2}    \left((\theta_{m',n'}
\phi_2)^{-1}(\mathfrak{M}_1 \mathfrak{N}) {\rm Nm}(\mathfrak{M}_1
\mathfrak{N})^{-\frac{z}{2}}\right ) \cdot
(\phi_2/\phi_1)(\mathfrak{N}) {\rm
Nm}(\mathfrak{N})^z.\end{eqnarray*}
\end{prop}

\begin{rem} Here the factor $\phi_2^{-1}(\mathfrak{M}_{1,v})$ at places $v
\in S$ stands for $\phi_{2,v}^{-r} (\pi_v)$ for the choice of
uniformizer $\pi_v$ in the definition of the newvector and
$\theta_{m',n'}(\mathfrak{N})$ for the product $\prod_{v \in T}
\theta_{m',n'}(\pi_v^{{\rm ord}_v \mathfrak{N}})$ for the
uniformizers $\pi_v$ chosen in the definition of ${\rm
Eis}^{\theta}(\Psi)$.
\end{rem}

\begin{proof}
 We start by unfolding the Eisenstein series $${\rm
Eis}^{\theta}(g,\Psi^{\rm new})= \sum_{\gamma \in B(\Q) \backslash
G(\Q)} \Psi^{\rm new, \theta}(\gamma g)$$ for ${\rm Re}(z) \gg 0$
and use analytic continuation to deduce the result for all $z$ for
which the integral converges.

Following \cite{Ko} \S4.5 we use a refinement of the Bruhat
decomposition choosing representatives for $B(\Q)\backslash G(\Q)$
according to the orbits of the $T(\Q)$-action:

$$G(\Q)=B(\Q)\begin{pmatrix}1&0 \\0 & 1 \end{pmatrix}
\cup B(\Q) w_0 \cup B(\Q) \begin{pmatrix}1&0 \\1 & 1 \end{pmatrix}
T_1(\Q),$$ where $T_1(\Q)=\left \{ \begin{pmatrix}1&0 \\0 & b
\end{pmatrix}: b \in F^* \right \}$.
If we decompose the integral according to this sum, the integral
over the first two summands vanishes, since $\omega_z(g_f
b_{\infty}, \phi, \Psi^{\rm new, \theta})=\Psi^{\rm new,
\theta}_f(g_f) \omega_{\infty} (b_{\infty})$ is zero along $H$
(here we factor (\ref{omega0}) as $\omega_z(g, \phi,
\Psi)=\omega_{\infty}(g_{\infty}) \cdot \Psi(g_f)$). We would like
to write the remaining term as

$$\int_{\A_F^*} \theta_{m',n'}(x) \Psi^{\rm new, \theta}(\begin{pmatrix}1&0 \\1 & 1 \end{pmatrix}
\begin{pmatrix}1&0 \\0 & x_f \end{pmatrix})
\cdot \left \langle Q_{m',n'},\omega_{\infty}(\begin{pmatrix}1&0
\\1 & 1 \end{pmatrix}\begin{pmatrix}1&0 \\0 & x_{\infty}
\end{pmatrix})(\frac{H}{2}) \right \rangle d^*x. $$
This step is justified if the latter integral converges
absolutely. Since the integrand decomposes by definition as a
product of local functions, the integral can be written as a
product of local integrals:

$$\prod_{v \nmid \infty} \int_{F_v^*} \theta_{m',n'}(x_v) \Psi^{\rm new, \theta}_v(\begin{pmatrix}1&0 \\1 & 1
\end{pmatrix} \begin{pmatrix}1&0 \\0 & x_v \end{pmatrix}) d^*
x_v$$
$$
\times \int_{F_{\infty}^*} \left \langle
Q_{m',n'},\omega_{\infty}(\begin{pmatrix}1&0
\\1 & 1 \end{pmatrix}\begin{pmatrix}1&0 \\0 & x_{\infty}
\end{pmatrix})(\frac{H}{2}) \right \rangle d^*x_{\infty}.$$

For each finite place $v$ we define integers $r=r_v,s=s_v$ by
$\mathfrak{P}_v^r \| \mathfrak{M}_1$ and $\mathfrak{P}_v^s \|
\mathfrak{M}_1 \mathfrak{M}_2$. We will treat the local integrals
according to the following cases:
\begin{enumerate}[(1)]
  \item $v$ finite place, $v \notin T$, both $\phi_i$ unramified, i.e., $r=s=0$
  \item $v$ finite place, $\phi_1$ ramified, $\phi_2$ unramfied,
  i.e., $r=s>0$
  \item $v$ finite place, $\phi_1$ unramified, $\phi_2$ ramified,
  i.e., $r=0$, $s>0$
  \item $v$ finite place, $r>0$ and $s-r>0$
  \item $v \in T$
  \item $v$ archimedean
\end{enumerate}

Before we start, we work out the Iwasawa decomposition of our
argument at the finite places:

$$\begin{pmatrix}1&0 \\1 & 1 \end{pmatrix} \begin{pmatrix}1&0 \\0 & x_v \end{pmatrix}=
\begin{pmatrix}1&0 \\1 & x_v \end{pmatrix}=
\begin{cases}
    \begin{pmatrix}x_v&1 \\0 & 1 \end{pmatrix} \begin{pmatrix}0&-1 \\1 &
    x_v\end{pmatrix}   & \text{ if } {\rm ord}_v(x_v) \geq0, \\

    \begin{pmatrix}1&0 \\0 & x_v \end{pmatrix}\begin{pmatrix}1&0 \\x_v^{-1} & 1
\end{pmatrix}  & \text{ if } {\rm ord}_v(x_v)<0.
\end{cases}
$$

We decompose $F_v^*$ into a disjoint union of $\pi_v^t
\mathcal{O}^*_{F_v}$ for $t \in \Z$.

 In case (1), the integrand over $\pi_v^t \mathcal{O}^*_v$  is
$$
  \begin{cases}
    (\phi_1 \theta_{m',n'})_v^t(\pi_v) |\pi_v|_v^{tz/2} & \text{ if } t \geq 0, \\
    (\phi_2 \theta_{m',n'})_v^t(\pi_v) |\pi_v|_v^{-tz/2} & \text{ if } t<0.
  \end{cases}
$$

 The integral therefore is given by two infinite sums, and since the
 infinity type of $\phi_1 \theta_{m',n'}$ is $z^{1+m-m'}
 \overline z^{-n'}$, and that of $\phi_2^{-1} \theta_{m',n'}^{-1}$ is
 $z^{1+m'} \overline z^{n'-n}$ it converges
 for $${\rm Re}(z)>n'-(m-m'+1) \text{ and } {\rm
 Re}(z)>(n-n')-(m'+1)$$ and the value is
\begin{eqnarray*}
&&\frac{1}{1-(\phi_{1,v}\theta_{m',n',v})(\pi_v) {\rm
Nm}(\mathfrak{P}_v)^{-z/2}} +
\frac{(\phi_{2,v}\theta_{m',n',v})^{-1}(\pi_v){\rm
Nm}(\mathfrak{P}_v)^{-z/2}}{1-
(\phi_{2,v}\theta_{m',n',v})^{-1}(\pi_v){\rm Nm}(\mathfrak{P}_v)^{-z/2}}\\
&=& \frac{ L_v(\frac{z}{2}, \phi_1 \theta_{m',n'})
L_v(\frac{z}{2}, (\phi_2 \theta_{m',n'})^{-1}) }{L_v(z,
\phi_1/\phi_2)}. \end{eqnarray*}

In case (2), the definition of the newvector $\Psi_v^{\rm new}$
shows that the integrand is non-zero only over $\pi_v^t
\mathcal{O}^*_v$ with $t\leq -r$. The integral therefore is given
by
\begin{eqnarray*}\lefteqn{\sum_{t\geq r} (\phi_{2,v}\theta_{m',n',v})^{-t}(\pi_v)
|\pi_v|_v^{tz/2}=} \\&=(\phi_{2,v}\theta_{m',n',v})^{-r}(\pi_v)
{\rm Nm}(\mathfrak{P}_v)^{-rz/2} \cdot L_v(\frac{z}{2}, (\phi_2
\theta_{m',n'})^{-1}).\end{eqnarray*}

For case (3) we only get a non-zero contribution for ${\rm
ord}_v(x_v) \geq 0$. Since for such $x_v$, $\begin{pmatrix}1&0 \\1
& x_v\end{pmatrix}=
\begin{pmatrix}x_v&* \\0 & 1 \end{pmatrix} \begin{pmatrix}1&0 \\1 & 1
\end{pmatrix}k$ with $k \in K^1(\mathfrak{P}_v^s)$, the
integral equals $$\sum_{t \geq 0}
(\phi_{1,v}\theta_{m',n',v})^t(\pi_v)
|\pi_v|_v^{tz/2}=L_v(\frac{z}{2}, \phi_1 \theta_{m',n'}).$$

In case (4), $\Psi_v$ is non-zero only on $\left [
\begin{pmatrix}1&0 \\\pi_v^r & 1 \end{pmatrix}\right ] \in B(F_v)
\backslash {\rm GL}_2(F_v) / K^1(\mathfrak{P}_v^s)$. This means we
have to have ${\rm ord}_v(x_v)=-r$ exactly. If $x_v= \epsilon
\pi_v^{-r}$ with $\epsilon \in \mathcal{O}_v^*$ we have

$$\begin{pmatrix}1&0 \\0 & x_v \end{pmatrix} \begin{pmatrix}1&0 \\x_v^{-1} & 1
\end{pmatrix}=\begin{pmatrix}1&0 \\0 & \epsilon \pi_v^{-r} \end{pmatrix}
\begin{pmatrix}1&0 \\\epsilon^{-1} \pi_v^{r} & 1 \end{pmatrix}=$$

$$\begin{pmatrix}1&0 \\0 & \epsilon \pi_v^{-r} \end{pmatrix}
 \cdot \begin{pmatrix}1&0 \\0 & \epsilon^{-1} \end{pmatrix} \begin{pmatrix}1&0 \\\pi_v^r & 1
\end{pmatrix} \begin{pmatrix}1&0 \\0 & \epsilon \end{pmatrix}=
\begin{pmatrix}1&0 \\0 & \pi_v^{-r} \end{pmatrix}
\begin{pmatrix}1&0 \\\pi_v^r & 1
\end{pmatrix} \begin{pmatrix}1&0 \\0 & \epsilon \end{pmatrix}.
$$

The integral therefore is given by
$$ \int_{\mathcal{O}_v^*} (\phi_{2,v}\theta_{m',n',v})^{-r} (\pi_v)
|\pi_v|_v^{rz/2} d^*\epsilon=(\phi_{2,v}\theta_{m',n',v})^{-r}
(\pi_v) {\rm Nm}(\mathfrak{P}_v)^{-rz/2}.$$

Case (5): For $v \in T$ Lemma \ref{thetatwist} implies that the
local factor is given by
$$\theta^{-1}_{m',n',v} (-1) (\phi_2/\phi_1)(\mathfrak{N}_v)
{\rm Nm}(\mathfrak{N}_v)^z L^{-1}_v(z, \phi_1/\phi_2) \int_{F_v^*}
\Psi^{\rm new}_{(\phi_1,\phi_2) |\alpha|^{z/2} \theta}
(\begin{pmatrix}1&0
\\1 & x_v
\end{pmatrix}) d^*x_v.
$$
Proceeding as in case (4) and taking the normalization of the
local measures into account we obtain
$$\theta_{m',n',v}^{-1} (-1)  L^{-1}_v(z, \phi_1/\phi_2) \cdot (\phi_1
\theta_{m',n'})_v^{-1}(\mathfrak{N}) {\rm
Nm}(\mathfrak{N}_v)^{z/2} \cdot
\#(\mathcal{O}_v/\mathfrak{N}_v)^*.$$

In Case (6)  the archimedean factor is
 $$ \frac{i}{8 \pi} \int_{\C^*} u^{m-m'+k} \overline u^{n-n'+\ell} \left \langle Q_{m',n'},
 \omega_{\infty}(\begin{pmatrix}
 1&0\\1&u \end{pmatrix}) (H) \right \rangle \frac{du \wedge d\overline u}{u
 \overline u}.$$
 Here we denote
$(\phi_{\infty} |\alpha|^{z/2}_{\infty}) (b_{\infty})
k_{\infty}^{-1}.\check{S}_+$ by
$\omega_{\infty}(b_{\infty}k_{\infty})$, so
$$\omega_{\infty}(b_{\infty}k_{\infty})(H)=((\phi_{1,\infty},
\phi_{2,\infty}) |\alpha|^{z/2}_{\infty}) (b_{\infty}) \check{S}_+
({\rm Ad}(k_{\infty})(H)) \otimes k^{-1}_{\infty}.Y^m \overline
X^n.$$

Our calculation of this factor essentially follows the one in
\cite{Ko} pp.111-113. One first obtains the Iwasawa decomposition
$$\begin{pmatrix} 1&0\\1&u
\end{pmatrix}=
\begin{pmatrix} \frac{u}{\sqrt{1+u \overline u}}& \frac{1}{\sqrt{1+u \overline u}}\\0&\sqrt{1+u \overline u} \end{pmatrix}
\begin{pmatrix} \frac{\overline u}{\sqrt{1+u \overline u}} & -\frac{1}{\sqrt{1+u \overline u}}\\ \frac{1}{\sqrt{1+u \overline u}}
 & \frac{u}{\sqrt{1+u \overline u}}\end{pmatrix}.$$
 We therefore get
$$(\phi_{\infty} |\alpha|^{z/2}_{\infty}) (\begin{pmatrix} \frac{u}{\sqrt{1+u \overline u}}&
 \frac{1}{\sqrt{1+u \overline u}}\\0&\sqrt{1+u \overline u} \end{pmatrix}) =
u^{1-k} \overline u^{-n-\ell} |u|_{\infty}^{z/2} \sqrt{1 + u
\overline u}^{n-m-2-2z}$$ Checking the action of $K_{\infty}$ on
the Lie algebra, we have
$$\check{S}_+(\begin{pmatrix} \frac{\overline u}{\sqrt{1+u
\overline u}} & -\frac{1}{\sqrt{1+u \overline u}}\\
\frac{1}{\sqrt{1+u \overline u}}
 & \frac{u}{\sqrt{1+u \overline u}}\end{pmatrix}. H)= \frac{2\overline u}{1+u \overline u}.$$
Lastly, we calculate
$$\left \langle Q_{m',n'}, k^{-1}_{\infty}.Y^m\overline X^n \right
\rangle= (-1)^{m-m'} \overline u^{m-m'+n'} \sqrt{1+u \overline
u}^{-m-n}.$$ Together this gives
$$\left\langle Q_{m',n'}, \omega_{\infty}(\begin{pmatrix}
 1&0\\1&u \end{pmatrix}) (H) \right \rangle=
 2 (-1)^{m-m'}\frac{(u \overline u)^{z/2 +1}\overline u^{m-m'-n+n'} }{u^k \overline u^{\ell}(1+u \overline u)^{z+2+m}}.$$
This gives rise to Beta-Function integrals, which converge for
$${\rm Re}(z/2)>-(m-m'+1), -(m'+1).$$ The archimedean integral
therefore contributes

\begin{eqnarray*}\lefteqn{(-1)^{m-m'}\frac{i}{4 \pi} \int_{\C^*} \frac{(u \overline
u)^{z/2 +1+m-m'}} {(1+u \overline u)^{z+2+m}} \frac{du \wedge
d\overline u}{u \overline u}}\\&=& \frac{(-1)^{m-m'}}{2}
\frac{\Gamma(z/2+m-m'+1)\Gamma(z/2+m'+1)}{\Gamma(z+m+2)}.\end{eqnarray*}

The preceding analysis shows that all the local integrals converge
absolutely for ${\rm Re}(z)\gg 0$ and that their product exists so
the integral over $\A_F^*$ converges absolutely.

To conclude the proof of the proposition by meromorphic
continuation it suffices to prove that for any $\xi \in G(\A_f)$
$$\int_{\sigma_{\xi} \otimes Q_{m',n'}} {\rm Eis}(\Psi^{\rm new, \theta})=
\int_0^{\infty} \theta_{m',n', \infty}(t) \left \langle
Q_{m',n'}, {\rm Eis}(\begin{pmatrix} 1 &0\\0& \xi t \end{pmatrix},
\Psi^{\rm new, \theta}_{\phi |\alpha|^{z/2}})(\frac{H}{2}) \right
\rangle \frac{dt}{t}$$ converges to a meromorphic function in $z$.
The following argument is adapted from \cite{CS} Proposition 3.5
and \cite{Wes} Proposition 2.4.

Recall that ${\rm Eis}(\Psi^{\rm new, \theta}_{\phi
|\alpha|^{z/2}})={\rm Eis}(\omega_z(\phi,\Psi^{\rm new,
\theta}_{\phi |\alpha|^{z/2}}))$. By picking out the
$\frac{\check{H}}{2}\otimes Q_{m',n'}^{\vee}$-component the
integrand equals $$t^{m-m'+n-n'+k+\ell} {\rm
Eis}(f(m',n',\Psi^{\rm new, \theta}_{\phi
|\alpha|^{z/2}})(\begin{pmatrix} 1 &0\\0& \xi t
\end{pmatrix}),$$ where $$f(m',n',\Psi^{\rm new,
\theta}_{\phi |\alpha|^{z/2}}) \in V_{\phi
|\alpha|^{z/2}}^{K_f^{\theta}} \otimes M^{m,n}_{\C}$$ is given by
$$(b_{\infty} k_{\infty} , g_f) \mapsto \phi_{\infty}(b_{\infty})
\tilde f(m',n',k_{\infty}) \Psi^{{\rm new}, \theta}_{\phi
|\alpha|^{z/2}}(g_f)$$ for a smooth function $\tilde f (m',n',
\cdot ):K_{\infty} \to \C^*$.

For $c>0$ let now $$I_c(z):= \int_{1/c}^c \theta_{m',n',
\infty}(t) {\rm Eis}(f(m',n',\Psi^{\rm new, \theta}_{\phi
|\alpha|^{z/2}})(\begin{pmatrix} 1 &0\\0& \xi t
\end{pmatrix}) \frac{dt}{t}.$$ Since the Eisenstein series has a meromorphic
continuation to all $z \in C$ this is a meromorphic function for
any $c>0$. It suffices therefore to show that $I_c(z)$ converges
locally uniformly for all $z$ as $c \rightarrow \infty$.

Put $E_z(g)={\rm Eis}(f(m',n',\Psi^{\rm new, \theta}_{\phi
|\alpha|^{z/2}})(g)$. One checks that the constant term ${\rm
res}(E_z) (g)$ vanishes for $g=\begin{pmatrix} 1 &0\\0& \xi t
\end{pmatrix}$ and $g=\begin{pmatrix} \xi &0\\0&  t
\end{pmatrix}w_0$.
It follows that $$I_c(z)=I_c^1(z) +I_c^2(z),$$ where
$$I_c^1(z)= \int_{1/c}^1 t^{m-m'+n-n'+k+\ell} \left( E_z(\begin{pmatrix} 1 &0\\0& \xi t
\end{pmatrix}) - {\rm res}(E_z) (\begin{pmatrix} 1 &0\\0& \xi t
\end{pmatrix})\right) \frac{dt}{t},$$
$$I_c^2(z)= \int_{1/c}^1 t^{m'+n'+k+\ell} \left( E_z - {\rm res}(E_z) \right) (\begin{pmatrix} 1 &0\\0& t
\end{pmatrix} \begin{pmatrix} \xi &0\\0& 1
\end{pmatrix} w_0)\frac{dt}{t}.$$
Standard growth estimates for automorphic forms on Siegel sets
(see \cite{Langl} Lemma 3.4, \cite{Schw} \S 1.10, \cite{HC} I
Lemma 10) imply that for any $g \in G(\A)$ and $r \in \Real$ there
exists a constant $C(g, r, z)>0$, locally uniform in $z$, such
that
$$|E_z(\begin{pmatrix} 1 &0\\0& t
\end{pmatrix} g) -{\rm res}(E_z)(\begin{pmatrix} 1 &0\\0& t
\end{pmatrix} g)| \leq C(g,r, z) t^r, \; 0<t \leq 1.$$ From this it follows that $I_c^1(z)$ and
$I_c^2(z)$ converge absolutely and locally uniformly for all $z$
as $c \rightarrow \infty$. The limits therefore define meromorphic
functions in $z$, as claimed above.

\end{proof}

\begin{cor} \label{torint}
For $n-1 < m'+n' < m+1$, $I(\phi, \theta, \Psi^{\rm new}, z)$
converges at $z=0$ and we get $I(\phi, \theta, \Psi^{\rm
twist}_{\phi_f}, 0)=$
\begin{eqnarray*} \lefteqn{
 \sum_{[\xi] \in \pi_0(K^{\theta}_f)} \theta_{m',n'}(\xi) \int_{\sigma_{\xi}
\otimes Q_{m',n'}} {\rm Eis}^{\theta}(\Psi^{\rm twist}_{\phi_f})=}
\\&=&\frac{L(0, \phi_1\theta_{m',n'}) L(0, \phi_2^{-1}
\theta^{-1}_{m',n'})} {L(0, \chi)} \cdot
\frac{\Gamma(m-m'+1)\Gamma(m'+1)}{\Gamma(m+2)} \cdot
C(\mathfrak{M}_1, S, \mathfrak{N}),\end{eqnarray*} where $$
C(\mathfrak{M}_1, S, \mathfrak{N})=
\frac{(-1)^{n-n'+k+\ell}}{2}(\theta_{m',n'}
\phi_2)^{-1}(\mathfrak{M}_1 \mathfrak{N}) \cdot
\chi^{-1}(\mathfrak{N}) \#(\mathcal{O}/\mathfrak{N})^*\cdot$$
$$\cdot \prod_{v \in S} (\mu_{2,v}^{-1} (-1)
\chi^{-1}(\mathfrak{P}_v^{r_v})).$$ \hfill \qedsymbol
\end{cor}

\section{Bounding the denominator}
After interpreting the toroidal integral as a cohomological
pairing we combine the calculation of Section 4 with results of
Hida and Finis to bound the denominator of the Eisenstein
cohomology class. For this we need the existence of certain Hecke
characters which we construct in Section \ref{s5.2}.
\subsection{Interpretation of the toroidal integral as evaluation
pairing}\label{s5.1} Let $F_{\phi, \theta}$ be the finite
extension of $F_{\phi}$ adjoining the values of the finite part of
$\theta_{m',n'}$ and the $L$-values $L^{\rm alg}(0, \phi_1
\theta_{m',n'})$ and $L^{\rm alg}(0, \phi_2^{-1}
\theta^{-1}_{m',n'})$ and denote its ring of integers by
$\mathcal{O}_{\phi,\theta}$. It follows from Proposition
\ref{ration} that $[{\rm Eis}^{\theta}(\Psi^{\rm twist}_{\phi_f})]
\in H^1(S_{K_f^{\theta}}, \widetilde{N^{\vee}_{F_{\phi,
\theta}}})$. Let
$$S_{K_f^{\theta}} \cong \bigoplus_{[\xi]\in \pi_0(K^{\theta}_f)}
\Gamma_{\xi}^{\theta} \backslash \h $$ for $\{[\xi]\}$ a system of
representatives of $\pi_0(K^{\theta}_f)$. Put
$$\partial_{\{0,\infty\},\xi}:=\Gamma_{\xi, B}^{\theta} \backslash e(B) \cup
\Gamma_{\xi, B^{w_0}}^{\theta} \backslash e(B^w) \subset
\Gamma_{\xi}^{\theta}\backslash \oh$$ and
$$\partial_{\{0,\infty\}}=\bigcup_{\xi} \partial_{\{0,\infty\},\xi}
\subset \partial \overline S_{K_f^{\theta}}.$$

 The relative cycles $\sigma_{\xi}\otimes Q_{m',n'}$ we described in \ref{3.2.1}
 give rise to classes in $$H_1(\Gamma_{\xi}^{\theta}
\backslash \oh, \partial_{\{0,\infty\},\xi}, j^*_{\xi} \widetilde
N_{\mathcal{O}_{\phi, \theta}}).$$

The following Lemma shows that ${\rm Eis}^{\theta}(\Psi^{\rm
twist}_{\phi_f})$ has vanishing constant terms at the $\infty$-
and $0$-cusps of each connected component. Since for any
automorphic form $f$ the function $f-f_P$ together with its
derivatives is fast decreasing at $P$  (see \cite{HC} I Lemma 10)
this implies that the cocycle gives rise to a differential form
fast decreasing at the $\infty$- and $0$-cusps of each connected
component. By Proposition \ref{Borel} the cocycle ${\rm
Eis}^{\theta}(\Psi^{\rm twist}_{\phi_f})$ therefore represents a
relative cohomology class in
$$\bigoplus_{[\xi] \in \pi_0(K^{\theta}_f)}
H^1(\Gamma_{\xi}^{\theta} \backslash \oh,
\partial_{\{0,\infty\},\xi}, j_{\xi}^*
\widetilde{N^{\vee}_{\C}}),$$ denoted by $[{\rm
Eis}^{\theta}(\Psi^{\rm twist}_{\phi_f})]_{\rm rel}$, mapping to
$[{\rm Eis}^{\theta}(\Psi^{\rm twist}_{\phi_f})] \in
H^1(S_{K_f^{\theta}},\widetilde{N^{\vee}_{F_{\phi,\theta}}})$.
This allows us to interpret the toroidal integral of the previous
section as sum of evaluation pairings for each connected component
so that the value of the integral provides a lower bound on the
denominator of $[{\rm Eis}^{\theta}(\Psi^{\rm
twist}_{\phi_f})]_{\rm rel}$ (see Section \ref{relativecoh} for
properties of the evaluation pairing).
\begin{lem}  \label{relativecocycle} We have $$[{\rm Eis}^{\theta}(\Psi^{\rm twist}_{\phi_f})]_{\rm rel} \in
H^1( S_{K_f^{\theta}} , \partial_{\{0,\infty\}},
\widetilde{N^{\vee}_{F_{\phi,\theta}}}).$$
\end{lem}
\begin{proof}
We claim that for $P=B$ and $P=B^{w_0}$  $${\rm
Eis}^{\theta}(\Psi^{\rm twist}_{\phi_f})_{P}(g_{\infty}
\begin{pmatrix}1 & 0\\ 0 & \xi
\end{pmatrix})=0$$ for all $g_{\infty} \in G_{\infty}$ and all $[\xi] \in \pi_0(K^{\theta}_f)$.
From the form of the constant term for ${\rm Eis}(\phi,
\Psi^0_{\phi_f})_B$ (see Proposition \ref{const}(a)) we deduce, by
interchanging the finite sums of the twists with the integral,
that
$${\rm Eis}^{\theta}(\Psi^{\rm twist}_{\phi_f})_B={\rm
Eis}^{\theta}(\omega_0(\phi, \Psi^{\rm
twist}_{\phi_f}))_B=\omega_0(\phi, (\Psi^{\rm
twist}_{\phi})^{\theta}) + c(\phi,0) \omega_0(w_0.\phi, (\Psi^{\rm
twist}_{w_0.\phi})^{\theta}).$$ We need to show that $(\Psi_*^{\rm
twist})^{\theta}$ vanishes on $\eta_f \begin{pmatrix} 1&0\\0&\xi
\end{pmatrix}$ for $\eta$ the identity matrix and $w_0$. Then vanishing
for $\eta$ equal to the identity matrix follows immediately from
$\sum_{x \in (\mathcal{O}_v/\mathfrak{N}_v)^*}
\theta^{-1}_{m',n',v}(x) =0$ for the finite order character
$\theta_{m',n',v}|_{\mathcal{O}_v^*}$ by definition of the
conductor $\mathfrak{N}_v$. For $\eta=w_0$ the vanishing follows
from our definition of the newvectors $\Psi^{\rm new}_*$, of which
$\Psi_*^{\rm twist}$ is a multiple, and from our choice of
$\mathfrak{N}$ coprime to the conductors of the characters
$\phi_1$ and $\phi_2$.

It remains to prove the rationality of $[{\rm
Eis}^{\theta}(\Psi^{\rm twist}_{\phi_f})]_{\rm rel}$. For this we
adapt an argument in \cite{CS} Lemma 5.2. Put
$$\omega=[{\rm Eis}^{\theta}(\Psi^{\rm twist}_{\phi_f})] \text{
and } \,\omega_{\rm rel}:=[{\rm Eis}^{\theta}(\Psi^{\rm
twist}_{\phi_f})]_{\rm rel}.$$ Let $$I_T=\{T_{\pi_v}
-\phi_2(\pi_v) - \phi_1(\pi_v) {\rm Nm}(\pi_v): v \notin S \cup T,
v\nmid \mathfrak{M} \}.$$ Then $I_T$ annihilates both $\omega$ and
$\omega_{\rm rel}$ (cf. Lemma \ref{Hecke}). Proposition
\ref{ration} implies that $\omega \in H^1( S_{K_f^{\theta}},
\widetilde{N^{\vee}_{F_{\phi,\theta}}})$. Thus, using a dimension
counting argument, $\omega$ is in the image of $$H^1(
S_{K_f^{\theta}} ,
\partial_{\{0,\infty\}},
\widetilde{N^{\vee}_{F_{\phi,\theta}}})[I_T] \to H^1(
S_{K_f^{\theta}}, \widetilde{N^{\vee}_{F_{\phi,\theta}}})[I_T],$$
where `$[I_T]$' denotes the subspaces annihilated by the elements
in $I_T$. Let $c$ be an element in the left hand side mapping to
$\omega$. Then $c-\omega_{\rm rel}\in H^1( S_{K_f^{\theta}} ,
\partial_{\{0,\infty\}},\widetilde{N^{\vee}_{\C}})[I_T]$ is in the
image of $H^0(\partial_{\{0,\infty\}},
\widetilde{N^{\vee}_{\C}})$. We recall the description of the
boundary cohomology as a $G(\A_f)$-module given by Harder:
$$H^0(\partial \overline S_{K_f^{\theta}}, \widetilde M_{\C})
\cong \bigoplus_{\mu:T(\Q) \backslash T(\A) \to \C^*}
V_{\mu,\C}^{K_f^{\theta}},$$ where in this case (degree 0) the sum
is over characters $\mu=(\mu_1,\mu_2)$ with infinity type (cf.
\cite{HaGL2} \S 3.5)
$$\mu_{1,\infty}(z)=z^{-m-k} \overline z^{-n-\ell} \text{ and } \,
\mu_{2,\infty}(z)=z^{-k} \overline z^{-\ell}.$$

By the Chebotarev density theorem we can find  an inert prime $q$
such that $q \equiv 1 \mod{\mathfrak{M}_1 \mathfrak{M}_2
\mathfrak{N}}$. We claim that $T_q=[K_f^{\theta} \begin{pmatrix} q
& 0\\ 0 &1\end{pmatrix}_q K_f^{\theta}]$ acts by a scalar factor
on $H^0(\partial \overline{S_{K_f^{\theta}}}, \widetilde M_{\C})$.
For this consider $\Psi \in V_{\mu,\C}^{K_f^{\theta}}$ for some
$\mu$ as above. As in Lemma \ref{Hecke} we get $$T_q
\Psi=(\mu_{2,q}(q) + q^2 \mu_{1,q}(q)) \Psi.$$ By our assumption
on $q$ it therefore acts by $q^{k+\ell} + q^2 \cdot
q^{m+n+k+\ell}$, independently of $\mu$. In particular, this also
describes the action of $T_q$ on $H^0(\partial_{\{0,\infty\}},
\widetilde{M_{\C}}) \subset H^0(\partial
\overline{S_{K_f^{\theta}}}, \widetilde M_{\C})$. Comparing this
with the fact that $T_q$ acts on $c-\omega_{\rm rel}$ by
$q^{m+k+\ell+1} + q^2 \cdot q^{n+k+\ell-1}$ shows that
$c-\omega_{\rm rel}=0$. This proves the rationality of
$\omega_{\rm rel}$.
\end{proof}

The next lemma shows that the denominator of the relative
Eisenstein cohomology class bounds the denominator of the original
Eisenstein cohomology class from below.
\begin{lem} \label{relden}
If $p> m\geq n$ then $$ \delta([\mathrm{Eis}^{\theta} (\Psi^{\rm
twist}_{\phi_f})]) \subseteq \delta([\mathrm{Eis}^{\theta}
(\Psi^{\rm twist}_{\phi_f})]_{\rm rel}).$$
\end{lem}
\begin{proof}
Put $\omega_{\rm rel}=[\mathrm{Eis}^{\theta} (\Psi^{\rm
twist}_{\phi_f})]_{\rm rel}$ and $\omega=[\mathrm{Eis}^{\theta}
(\Psi^{\rm twist}_{\phi_f})]$. Suppose $a \in
\mathcal{O}_{\phi,\theta}$ is such that $$a \cdot \omega \in
H^1(S_{K_f^{\theta}} ,
\widetilde{N^{\vee}_{\mathcal{O}_{\phi,\theta}}})_{\rm free}
\text{ but } \, a \cdot \omega_{\rm rel} \notin
H^1(S_{K_f^{\theta}},
\partial_{\{0,\infty\}},\widetilde{N^{\vee}_{\mathcal{O}_{\phi,\theta}}})_{\rm free}.$$
Let $\lambda$ be a uniformizer of $\mathcal{O}_{\phi,\theta}$ and
$m \geq 1$ the smallest integer such that $\lambda^m a \omega_{\rm
rel}$ is in the image of an element, say $c$, of
$H^1(S_{K_f^{\theta}},\partial_{\{0,\infty\}},\widetilde{N^{\vee}_{\mathcal{O}_{\phi,\theta}}})$.
Then the image $\overline c$ of $c$ in $H^1(S_{K_f^{\theta}},
\partial_{\{0,\infty\}},\widetilde{N^{\vee}_{k}})$ (where
$k=\mathcal{O}_{\phi,\theta}/\lambda$) is nonzero, but its image
in $H^1(S_{K_f^{\theta}},\widetilde{N^{\vee}_{k}})$ is zero.
Therefore $\overline c$ is in the image of
$H^0(\partial_{\{0,\infty\}}, \widetilde{N^{\vee}_{k}})$. By
Nakayama's Lemma it even has to be in the image of
$$H^0(\partial_{\{0,\infty\}}, \widetilde{N^{\vee}_{k}})/
H^0(\partial_{\{0,\infty\}},\widetilde{N^{\vee}_{\mathcal{O}_{\phi,\theta}}})\otimes
k.$$ Note that this quotient is isomorphic to the
$\lambda$-torsion of
$H^1(\partial_{\{0,\infty\}},\widetilde{N^{\vee}_{\mathcal{O}_{\phi,\theta}}})$.
Under our assumptions $p$ does not divide the level $K_f$ (i.e.,
$K_v={\rm GL}_2(\mathcal{O}_v)$ for $v \mid p$) so Proposition
2.4.1 (ii) of \cite{U95} shows that if $p>{\rm max}\{n,3\}$ then
$H^1(\partial_{\{0,\infty\}},\widetilde{N^{\vee}_{\mathcal{O}_{\phi,\theta}}})$
is torsion-free (the argument in \cite{U95} extends to our more
general coefficient system). This shows that $\overline c=0$, in
contradiction to our assumption, so $a$ cannot exist, proving our
Lemma.
\end{proof}

\subsection{Construction of special Hecke characters} \label{s5.2}
Recall that for a Hecke character $\lambda: F^*\backslash \A_F^*
\to \C^*$ we defined $\lambda^*(x)=\lambda^{-1} (\overline x)
|x|$. Following constructions by Greenberg \cite{G85}, Rohrlich
\cite{Ro}, and Yang \cite{Yang} we prove the following:

\begin{lem} \label{greenchar}

  \emph{(a)} For $F= \Q(\sqrt{-1})$ or $\Q(\sqrt{-3})$ there exists a Hecke
  character $\mu^{(1,0)}$ of infinity type $z$ with conductor   $2 \mathcal{D}$ such that
  $(\mu^{(1,0)})^*=\mu^{(1,0)}$.

\emph{(b)} If $F \neq \Q(\sqrt{-1}),\Q(\sqrt{-3})$ then for any
$k>0$
  there exists a Hecke character
$\mu^{(k,1-k)}$ of infinity type $z^k \overline z^{1-k}$ such that
$(\mu^{(k,1-k)})^*=\mu^{(k,1-k)}$ whose conductor is given by $$
\begin{cases}
    \mathcal{D} & \text{ if } d_F \text{ odd}, \\
    2 \mathcal{D} & \text{ if } d_F \text{ even}.
  \end{cases}$$

\end{lem}

\begin{proof}
For $F=\Q(\sqrt{-1})$ and $\Q(\sqrt{-3})$ one can take the inverse
of the Gr\"ossen-characters associated to the elliptic curves
$y^2=x^3+x$ (conductor 64) or $y^2=x^3+1$ (conductor 36),
respectively (for curves with minimal conductor divisible only by
ramified primes see \cite{G85} Lemma p.81). For $F \neq
\Q(\sqrt{-1}),\Q(\sqrt{-3})$ we note that Greenberg's construction
can be extended to $k\geq 1$: Let $p_1, p_2, \ldots, p_t$ be the
rational primes dividing the discriminant $d_F$ and let
$\mathfrak{p}_1, \mathfrak{p}_2, \ldots, \mathfrak{p}_t$ be the
corresponding primes of $F$. Since ${\rm
Nm}(\mathcal{O}^*_{\mathfrak{p}_i})$ is of index 2 in $\Z_{p_i}^*$
one can define a character  of order 2 on $\Z_{p_i}^*$
 with kernel containing ${\rm Nm}(\mathcal{O}^*_{\mathfrak{p}_i})$.
Via the embedding of $\Z_{p_i}^* \hookrightarrow
\mathcal{O}^*_{\mathfrak{p}_i}$ this character can be extended to
a character $\Psi_i $ of $\mathcal{O}^*_{\mathfrak{p}_i}$ having
finite order (can choose order 2 unless $p_i=2$ and $4 \| d_F$).
We can therefore define a continuous homomorphism $\Psi: \C^*
\cdot \prod_v \mathcal{O}_v^* \to \C^*$ so that $\Psi(z)=z^k
\overline z^{1-k}$ for $z \in \C^*$,
$\Psi|_{\mathcal{O}_{\mathfrak{p}_i}^*}=\Psi_i$ for $1 \leq i \leq
t$, and $\Psi$ is trivial on the other local units. Since $-1$ is
the only non-trivial unit and $\Psi(-1)=1$ we can define $\Psi$ to
be trivial on $F^*$. This character $\Psi$ can now be extended to
a Hecke character $\mu^{(k,1-k)}$ on $\A_F^*$.

We check that $(\mu^{(k,1-k)})^*=\mu^{(k,1-k)}$ by showing that
$\mu^{(k,1-k)}|_{\A^*}= \omega_{F/\Q} |\cdot|_{\A^*}$ for
$\omega_{F/\Q}$ the quadratic character associated to $F/\Q$ (see
\cite{G85} for a different proof): Clearly $(\mu^{(k,1-k)} | \cdot
|_{\A}^{-1})(t)=1$ for $t \in \Real^*_{>0}$ and $t \in \Q^*$, but
$$(\mu^{(k,1-k)} | \cdot |_{\A}^{-1})(-1)=-1.$$ By construction it
is also trivial on ${\rm Nm}(\A_F^*)$.

At odd primes the conductor is clearly of index 1. For the
calculation of the conductor at the place dividing 2 (and the
existence of characters with conductors as claimed) see Rohrlich
\cite{Ro} ($8 \mid d_F$) and Yang \cite{Yang} ($4 \| d_F$). Note
that we do not take one of Yang's characters with minimal
conductor but one with index 4 at the prime dividing 2.
\end{proof}

\begin{remno} \label{lambdastar}
\begin{enumerate}
\item We note that any algebraic Hecke character $\lambda$
satisfying $\lambda^*=\lambda$ is of the form $\mu^{(k,1-k)} \cdot
\vartheta$ for a finite order anticyclotomic character $\vartheta$
(i.e., such that $\vartheta^c=\overline \vartheta=\vartheta^{-1}$)
and that they satisfy $\lambda|_{\A^*}=\omega_{F/\Q} |
\cdot|_{\A^*}$ with $\omega_{F/\Q}$ the quadratic character of
$\Q^*\backslash \A^*$ associated to $F/\Q$.

\item More generally, for unitary Hecke characters $\lambda$
satisfying $\lambda^c=\overline \lambda$ we have
$$\lambda|_{\A^*}=
  \begin{cases}
    1 & \text{ if } \lambda_{\infty}(-1)=1, \\
    \omega_{F/\Q} & \text{ if } \lambda_{\infty}(-1)=-1.
  \end{cases}
$$
\end{enumerate}
\end{remno}

In addition,  we note the existence of the following character
(cf. \cite{Ti} Lemme 2.5, \cite{dS} Lemma II.1.4(ii)):
\begin{lem} \label{minram}
Let $q \geq 5$ be a rational prime and $\mathfrak{q}$ a prime of
$F$ dividing $q$. Then there exists a Hecke character with
conductor $\mathfrak{q}$ of infinity type $z$.
\end{lem}
\begin{proof}
Since $q \geq 5$, $\mathfrak{q}$ separates the roots of unity and
so the character is well-defined on $F^* \cdot \C^*
U(\mathfrak{q})$, where $U(\mathfrak{q}):=\{x \in
\hat{\mathcal{O}}^* | x \equiv 1 \mod{\mathfrak{q}
\hat{\mathcal{O}}} \}.$ Since the ray class group $ F^* \backslash
\A_{F,f}^* / U(\mathfrak{q})$ is finite we can  extend trivially
to a continuous character on $\A_F^*$.
\end{proof}

\subsection{Bounding the denominator} \label{3.5}
Because of Lemma \ref{relden} we now assume in addition that
$p>m$. We are interested in bounding
$$\delta([\mathrm{{\rm Eis}} (\Psi^0_{\phi_f})])=
\{a \in \mathcal{O}_{\phi}: a\cdot [\mathrm{{\rm Eis}} (
\Psi^0_{\phi_f})] \in H^1( S_{K_f^{S}}, \widetilde{
(N_{\mathcal{O}_{\phi}})^{\vee}})_{\rm free}  \} .$$ Observe that
$$\delta([\mathrm{{\rm Eis}} (\Psi^0_{\phi_f})]) \subseteq \delta([\mathrm{{\rm Eis}} (
\Psi^{\rm twist}_{\phi_f})])\subseteq \mathcal{O}_{\phi},$$
$$\delta([\mathrm{Eis} ( \Psi^{\rm twist}_{\phi_f})]) \mathcal{O}_{\phi, \theta}
\subseteq \delta([\mathrm{Eis}^{\theta} (\Psi^{\rm
twist}_{\phi_f})]) \subset \mathcal{O}_{\phi, \theta},$$ and (by
Lemma \ref{relden})
$$ \delta([\mathrm{Eis}^{\theta} (\Psi^{\rm
twist}_{\phi_f})]) \subseteq \delta([\mathrm{Eis}^{\theta}
(\Psi^{\rm twist}_{\phi_f})]_{\rm rel}).$$
 In Section \ref{s5.1} we showed that the
toroidal integral
$$I(\phi, \theta, \Psi^{\rm twist}_{\phi_f}, 0 )=\sum_{[\xi] \in
\pi_0(K^{\theta}_f)} \theta_{m',n'}(\xi) \int_{\sigma_{\xi}\otimes
Q_{m',n'}} {\rm Eis}^{\theta}(\Psi^{\rm twist}_{\phi_f}) $$ gives
the value of a sum of evaluation pairings between relative
cohomology and homology. The functoriality of these pairings
implies that the denominator  of $[\mathrm{Eis}^{\theta}
(\Psi^{\rm twist}_{\phi_f})]_{\rm rel}$ is bounded below by the
denominator of the integral. From Corollary \ref{torint} we deduce
that
$$I(\phi, \theta,
\Psi^{\rm twist}_{\phi_f}, 0)= \frac{L^{\rm alg}(0,\phi_1
\theta_{m',n'}) L^{\rm alg}(0, (\phi_2 \theta
_{m',n'})^{-1})}{L^{\rm alg}(0,\phi_1/\phi_2)} \cdot
C(\mathfrak{M}_1, S, \mathfrak{N}).$$

Since the conductors of $\phi_i$ and $\theta_{m',n'}$ and
$\#(\mathcal{O}/\mathfrak{N})^*$ are coprime to $(p)$ one checks
using Lemma \ref{charinteg} that $C(\mathfrak{M}_1, S,
\mathfrak{N}) \in \mathcal{O}_{\phi, \theta}^*$. This shows that
$\delta([\mathrm{Eis}^{\theta} (\Psi^{\rm twist}_{\phi_f})]_{\rm
rel})$ is contained in the (possibly fractional) ideal
$$\left ( \frac{L^{\rm alg} (0,\phi_1/\phi_2)}{L^{\rm alg}(0,\phi_1 \theta_{m',n'}) L^{\rm alg}(0, (\phi_2
\theta)_{m',n'}^{-1})} \right )  \mathcal{O}_{\phi, \theta}.$$
\begin{prop} \label{principle}
If there exists a Hecke character $\theta_{m',n'}$ as in
(\ref{theta}) such that $L^{\rm alg}(0,\phi_1 \theta_{m',n'})$ and
$L^{\rm alg}(0, (\phi_2 \theta_{m',n'})^{-1})$ lie in
$\mathcal{O}^*_{\phi,\theta}$ then
$$\delta([\mathrm{Eis} (\Psi^0_{\phi_f})]) \subseteq L^{\mathrm{alg}}(0, \chi)
\mathcal{O}_{\phi}$$ for $\chi=\phi_1/\phi_2$.
\end{prop}

We have at our disposal two results on the non-vanishing modulo
$p$ of Hecke $L$-values  as the Hecke character varies in an
anticyclotomic $\Z_q$-extension for $q \neq p$:
\begin{thm}[(Finis \cite{Fi2} Thm. 1.1)] \label{Finis}
Let $q$ be an odd prime split in $F$, distinct from $p$.
 Consider Hecke characters $\lambda$ of
infinity type $\lambda_{\infty}(z)=z^a \overline z^{1-a}$ for a
fixed positive integer $a$ with $\lambda^*=\lambda$, conductor
dividing $dd_Fq^{\infty}$ for some fixed $d$ coprime to $(p)$,
global root number $W(\lambda)=1$, and such that no inert primes
congruent to $-1  \mod{p}$ divide the conductor of $\lambda$ with
multiplicity one. If $a>1$ then assume $p$ splits in $F$. Then for
all but finitely many such Hecke characters
$$L^{\rm alg}(0,\lambda)  \text{ is a } p-\text{adic unit}.$$
\end{thm} Hida has proved a
similar result:
\begin{thm}[(\cite{Hi04b} Theorem 4.3)] \label{Hida}
Assume $p$ splits in $F$. Fix a character $\lambda$ of split
conductor (i.e., such that the conductor is a product of primes
split in $F/\Q$) coprime to $p$ with infinity type
$\lambda_{\infty}(z)=z^a \left(\frac{z}{\overline z}\right)^b$ for
$a>0$ and $b \geq 0$. Let $q$ be a split prime distinct from $p$
and coprime to the conductor of $\lambda$. Then
$$ L^{\rm alg}(\lambda \vartheta,0) \text{ is a } p-\text{adic unit}$$
for all but finitely many finite-order anticyclotomic characters
$\vartheta$ of $q$-power conductor.
\end{thm}

\begin{rem}
We quoted above the cases of Finis' Theorem when all but finitely
many $L$-values in the anticyclotomic tower are $p$-adic units; in
general this is not true, see \cite{Fi2} for the full statement.
Finis also allows ramification at $p$. Hida's Theorem is actually
valid for general CM-fields and also treats the case of non-split
$q$.
\end{rem}

We can show then, for example, the following:
\begin{thm} \label{thm01}
If $p>m$ is split in $F$ and both $\phi_i$ have split conductor
coprime to $(p)$, then
$$\delta([\mathrm{Eis} (\Psi^0_{\phi_f})]) \subseteq L^{\mathrm{alg}}(0, \chi)
\mathcal{O}_{\phi}.$$
\end{thm}

\begin{proof}
By Lemma \ref{minram} we can always find a character
$\theta_{m',n'}$ of the correct infinity type with split conductor
$\mathfrak{N}$ coprime to $(p d_F)$ and the conductors of the
$\phi_i$ such that $\#(\mathcal{O}/\mathfrak{N})^*$ is also
coprime to $(p)$. Applying Theorem \ref{Hida} for both $\phi_1
\theta_{m',n'}$ and $\phi_2^{-1} \theta_{m',n'}^{-1}$  there
exists a split prime $q \neq p$ coprime to the conductors of the
$\phi_i$ with $q \not \equiv 1 \mod{p}$ and a finite order
anticyclotomic character $\vartheta$ of $q$-power conductor such
that $L^{\mathrm{alg}*}(0, \phi_1 \theta_{m',n'} \vartheta )$ and
$L^{\mathrm{alg}*}(0, (\phi_2 \theta_{m',n'} \vartheta)^{-1})$
both lie in $\mathcal{O}_{\phi, \theta \vartheta}^*$ and we can
apply Proposition \ref{principle} for this modified character
$\theta'_{m',n'}=\theta_{m',n'} \vartheta$.
\end{proof}

\begin{rem}
This is where our restriction to $m \geq n$ is needed so that the
infinity types of $\phi_1 \theta_{m',n'}$ and $\phi_2^{-1}
\theta_{m',n'}^{-1}$ satisfy the condition of Theorem \ref{Hida}.
By using the $p$-adic functional equation it might be possible to
extend Hida's result to $a \leq 1$ and $b \geq 1-a$, which would
remove this condition.
\end{rem}

For finding congruences between the Eisenstein cohomology class,
multiplied by its denominator, and a cuspidal cohomology class, as
described in the introduction, we are interested in the case when
the restriction of the Eisenstein class to the boundary is
integral. As described in Proposition \ref{resinteg} we know that
this is the case when $m=n$ and $\chi^c=\overline \chi$. The
following theorem shows that in this situation there exists (under
some conditions on the conductor of $\chi$) an Eisenstein
cohomology class with $L^{\rm alg}(0,\chi)$ as lower bound on the
denominator. Note that for $m=n>0$ the results of Hida and Finis
are applicable only for primes $p$ split in $F$. Recall the
definition of the Gauss sum $\tau(\tilde \chi)$ from Section
\ref{Hcharacters}.

\begin{thm} \label{thm02}
Let $\chi$ be a Hecke character  of infinity type $z^{m+2}
\overline z^{-n}$  for $m \geq n \in \N_{\geq 0}$ with conductor
$\mathfrak{M}$ coprime to $(p)$. Assume $p>m$ and in addition that
either

(i) $p$ splits in $F$ and $\chi$ has split conductor

or

 (ii) $m=n$, (if $m>0$ then also assume that $p$ is split), $\chi^c=\overline \chi$, no
ramified primes (or 2
  if $F=\Q(\sqrt{-3})$)
  divide $\mathfrak{M}$ and no inert primes
congruent to $-1 \mod{p}$ divide $\mathfrak{M}$ with multiplicity
one, and
$$\omega_{F/\Q}(\mathfrak{M}) \frac{\tau(\tilde \chi)}{\sqrt{{\rm Nm}(\mathfrak{M})}}=1.$$

Then there exists a character $\phi=(\phi_1, \phi_2)$ with
$\chi=\phi_1/\phi_2$ such that the conductor of $\phi_1$ is
coprime to $(p) \mathfrak{M}$ and
$$\delta([\mathrm{Eis} (\Psi^0_{(\phi_1, \phi_2)_f})]) \subseteq
(L^{\mathrm{alg}}(0, \chi)).$$

\end{thm}

\begin{proof}
Part (i) follows directly from Theorem \ref{thm01} and Lemma
\ref{minram}. For (ii) we choose $m'=m$ and $n'=0$ (so that the
toroidal integral converges) and $k=0$, $\ell=-m$. This means that
$\theta_{m',n'}$ has to be a finite order character and $\phi_1$
should have infinity type $z$. For suitable $\phi_1$ and $\phi_2$
we want to apply Theorem \ref{Finis} to find a finite order
anticyclotomic $\theta_{m',n'}$ of $q$-power conductor, $q \neq p$
split prime coprime to the conductors of the $\phi_i$ and $q \not
\equiv 1 \mod{p}$, such that both $L^{\rm alg}(0, \phi_1
\theta_{m',n'})$ and $L^{\rm alg}(0,(\phi_2 \theta_{m',n'})^{-1})$
lie in $\mathcal{O}_{\phi \theta}^*$ and such that
\begin{equation} \label{root} W(\phi_1 \theta_{m',n'})=W((\phi_2
\theta_{m',n'})^{-1})=1.\end{equation} The characters $\phi_i$
must satisfy the conditions on the conductor imposed in Theorem
\ref{Finis} and $\phi_1^*=\phi_1$, $\phi_2^{-1}=(\phi_2^{-1})^*$.

Furthermore, (\ref{root}) imposes a condition on the root numbers
of the $\phi_i$ as we will now show: Let $\lambda$ be any Hecke
character satisfying $\lambda^*=\lambda$ with conductor
$\frak{f}_{\lambda}$ and $\vartheta$ a finite order anticyclotomic
character with conductor $Q^n$ for $Q \in \Z$ prime, $Q \neq 2$
and coprime to $\mathfrak{f}_{\lambda}$. Since
$\mathfrak{f}_{\lambda}=\overline{\mathfrak{f}}_{\lambda}$ we get
$\vartheta(\mathfrak{f}_{\lambda})= \pm 1$, but by assumption
$\vartheta$ has only $Q$-power roots of unity as values, so
$\vartheta(\mathfrak{f}_{\lambda})=1$. Also it is known that
$W(\vartheta)=1$ (see, for example, \cite{G83} p. 247 and
\cite{FQ}). By Remark \ref{lambdastar} we know $\tilde
\lambda(Q^n)=\omega_{F/\Q} (Q^n)$ for $\tilde
\lambda=\lambda/|\lambda|$. Proposition \ref{rootprod} therefore
shows that
\begin{equation} \label{rootprod2} W(\lambda \vartheta)=W(\lambda)
W(\vartheta) \tilde \lambda(Q^n)
\vartheta(\mathfrak{f}_{\lambda})=W(\lambda) \omega_{F/\Q}(Q^n).
\end{equation}
This implies that we need  $W(\phi_1)=W(\phi_2^{-1})=1$ to be able
to satisfy (\ref{root}) because we are considering $Q=q$ split.

We now define $\phi_1$: By possibly twisting $\mu^{(1,0)}$ from
Lemma \ref{greenchar} by a finite order anticyclotomic character
$\vartheta$ with suitable inert conductor we can always ensure by
(\ref{rootprod2}) that the resulting character, which we take as
$\phi_1$, satisfies $\phi_1^*=\phi_1$, $\phi_{1,\infty}(z)=z$,
$W(\phi_1)=1$, and ${\rm cond}(\phi_1)=r \mathcal{D}$, for $r \in
\Z$ coprime to $p \, \mathfrak{M}$ and such that no inert prime
$\equiv -1 \mod{p}$ divides $r$ with multiplicity one.

One checks that under our assumptions $\chi^c=\overline \chi$ and
$m=n$
 the character $\phi_2^{-1}=\chi/\phi_1$ satisfies
$(\phi_2^{-1})^*=\phi_2^{-1}$. From the definition in Section
\ref{Hcharacters} we deduce that $W(\chi)=- \frac{\tau(\tilde
\chi)}{\sqrt{{\rm Nm}(\mathfrak{M})}} \tilde
\chi(\mathcal{D}^{-1})$. Now applying Proposition \ref{rootprod}
we calculate that
$$W(\phi_2^{-1})=W(\phi_1^{-1} \chi)=-W(\phi_1^{-1}) W(\chi)
\omega_{F/\Q}(\mathfrak{M}) \tilde \chi (r \mathcal{D})
\overset{\text{assumption}}{=}W(\phi_1^{-1})=W(\phi_1),$$ as
desired. Here we use again Remark \ref{lambdastar} ($\tilde
\phi_1|_{\A^*}=\omega_{F/\Q}$ and $\tilde \chi|_{\A^*}\equiv 1$),
and the last equality holds because $\phi_1^c=\overline \phi_1$.
By Theorem \ref{Finis} there exists now some finite order
character $\theta_{m',n'}$ such that $L^{\rm alg}(0, \phi_1
\theta_{m',n'})$ and $L^{\rm alg}(0,(\phi_2 \theta_{m',n'})^{-1})$
are simultaneously $p$-adic units.
\end{proof}

\begin{rem} The condition
$\omega_{F/\Q}(\mathfrak{M}) \frac{\tau(\tilde \chi)}{\sqrt{{\rm
Nm}(\mathfrak{M})}}=1$ is satisfied, for example, by everywhere
unramified characters, so the theorem holds for any split or inert
prime $p$ and unramified $\chi$ with infinity type $z^2$.
\end{rem}

\bibliographystyle{amsplain}
\bibliography{biblio3}
\end{document}